%% file: AAHRah.tex
\documentclass[11pt]{article}
\usepackage[novbox]{pdfsync}

\usepackage[margin=2cm]{geometry}
\usepackage{mathrsfs,float,pgfplots,
bbm,graphicx,subcaption,listings,xcolor,tcolorbox,bookmark,enumitem,algorithmicx,algorithm,
algpseudocode,booktabs,array,multirow}
\usepackage[version=4]{mhchem}
\usepackage{tikz-cd}
\usepackage{lineno}
\tikzcdset{every label/.append style = {font = \small}}
\usepackage{tikz}
\usetikzlibrary{graphs,decorations.pathmorphing,decorations.markings,arrows.meta,
automata,positioning}
\usepackage{nomencl}
\makenomenclature

\usepackage{hyperref}
\usepackage[toc,sort=use]{glossaries}
\makeglossaries
\providecommand{\email}[1]{\texttt{#1}}

\input{AAdefN}

\renewcommand{\printboth}{%
  \printglossary[style=unifiedstyle, title={0 \hspace{.22em} Definitions, Theorems, Lemmas, Remarks, Examples  and Open Problems}]%
}

\def\PV{Perron-Volterra candidate Lyapunov function}
\title{A  Perron-Volterra Lyapunov function for  mathematical epidemiology reaction networks (MERN)   with non-interacting  rank-one strains}

\author{
Rim Adenane\thanks{Laboratoire d'Analyse, G\'eom\'etrie et Applications, D\'epartement des Math\'ematiques, Universit\'e Ibn-Tofail, K\'enitra, 14000, Morocco. \texttt{rim.adenane@uit.ac.ma}},
Florin Avram\thanks{Department de Math\'emathiques, Universit\'e de Pau, France. (\email{avramf3@gmail.com})
}, Miruna Beldiman\thanks{},
Andrei-Dan Halanay\thanks{Department of Mathematics and Computer Science, University of Bucharest, Bucharest, RO-010014, Romania. (\email{halanay@fmi.unibuc.ro})
}
}

\pgfplotsset{compat=1.18}

\begin{document}

\maketitle
\begin{abstract}
Persistence, coexistence,  competitive exclusion and global asymptotic stability
(GAS) are important related problems in mathematical epidemiology (ME), which are often studied through model-specific analyses.   Below, we provide a general approach to the last problem, based on constructing Perron-Volterra Lyapunov functions $L_p$, \eqr{Hp}, which apply to a general class of mathematical epidemiology models. More precisely, these Lyapunov functions yield  competitive-exclusion GAS partitions of the parameter space, for some bilinear m-strain models  with irreducible rank-one infection blocks, which do not interact, in the sense that the next  generation matrix has a block-diagonal  structure. There are at most  $m+1$ generic partition sets, corresponding to exactly one dominant strain, or to none, in which case the \DFE\ (DFE) is GAS -- see Theorem \ref{t:blo}. Coexistence equilibria with two or more strains may occur only on non-generic tie surfaces where  corresponding reproduction numbers are equal.
A second general  competitive-exclusion GAS partition -- see Theorem \ref{t:GAS} is provided for two-strain models with increasing concave incidence and scalar, uncorrelated strain blocks.  Here, equilibria supported on all $2^m$ boundary faces are possible.
Our approach rests on several pillars:
(i)   the fact that forward‑invariant boundary faces are determined by  siphon sets  of variables (a \CRN\ and Petri nets concept), with the disease‑free face being  the intersection of all minimal
 siphons, and the fact that Jacobians on siphon faces have a triangular block form, cf.
 Avram, Adenane, and Halanay(2026). (ii) The recently  established fact in Avram, Adenane, Halanay, Horvath and Khong (2026) that a transversal Jacobian block on a siphon face is Metzler, which puts under spotlight the roles of its Perron eigenvectors.
 (iii) \NGMs\ (NGMs), a \ME\ concept: on siphon faces, each regular
splitting of the transversal Jacobian allows defining a  NGM, and invasibility may be determined by comparing its spectral radium to $ 1$; this is particular useful when the NGM has rank one, like throughout this paper.   (iv) Boundary transcritical invasion relays (an ecology concept):  a bifurcation theorem linking eigenvalue crossing to the emergence of
a positive branch on an adjacent face.
(v) A  \PV\ combining Volterra entropy terms for the resident variables on a face with  Perron-weighted linear functionals for the invading blocks, whose weights are left Perron eigenvectors of corresponding NGMs.
We have provided an algorithmic implementation in the Mathematica
package \textsc{EpidCRN}
(\url{https://github.com/florinav/EpidCRNmodels}). The method organizes
the dynamics via the minimal siphon lattice, recursively computes Perron
eigenvectors of transversal Jacobians, and constructs candidate Lyapunov
functions for all equilibria, producing a parameter partition into
regions with a unique locally stable equilibrium. In the two cases mentioned above,
the candidate Lyapunov
functions are indeed  Lyapunov
functions.
\end{abstract}


\tableofcontents
\printboth
\input{intR1}
\input{mt}

\input{Fr}

\input{ct}
\input{wh}
\input{GASt}

\input{not}

\input{GASpr}
\input{alAH}

\input{remOP}

\input{R2}

\input{RahE2AH}
\input{exa}
\input{conc}

\bibliographystyle{amsalpha}
\bibliography{ref}
\glsaddall  

\end{document}

%% file: AAdefN.tex
\RequirePackage{etoolbox}
\usepackage{makecell}
\usepackage{amssymb, amsfonts,mathtools,bbm,eurosym,graphicx,latexsym, amsmath,amsthm,times,nccmath,url}
\usepackage{listings,verbatim,cases,hyperref,tabularx}
\usepackage{graphicx,float,fancybox}
\usepackage{epsfig,epstopdf,tikz}
\usepackage[sans]{dsfont}
\usetikzlibrary{positioning,arrows,arrows.meta,calc}
\lstdefinestyle{mathematica}{
 language=Mathematica,
 basicstyle=\small\ttfamily,
 backgroundcolor=\color{gray!10},
 frame=single,
 breaklines=true,
 commentstyle=\color{green!60!black},
 keywordstyle=\color{blue},
 stringstyle=\color{red},
 showstringspaces=false,
 tabsize=2
}
\usepackage{tcolorbox}
\newtcolorbox{keyresult}[1][]{
 colback=blue!5!white,
 colframe=blue!75!black,
 title=#1,
 fonttitle=\bfseries
}
\newtcolorbox{examplebox}[1][]{
 colback=green!5!white,
 colframe=green!75!black,
 title=#1,
 fonttitle=\bfseries
}

\makeatletter
\@ifundefined{theorem}   {\newtheorem{theorem}   {Theorem}}   {}
\@ifundefined{definition}{\newtheorem{definition}{Definition}} {}
\@ifundefined{open}      {\newtheorem{open}      {Problem}}   {}
\@ifundefined{example}   {\newtheorem{example}   {Example}}   {}
\makeatother

\newcounter{tabcounter}
\newcounter{defcounter}
\newcounter{opcounter}
\newcounter{thmcounter}
\newcounter{excounter}
\newcounter{ascounter}
\newcounter{cjcounter}
\newcounter{lemcounter}
\newcounter{remcounter}
\newcounter{figcounter}

\newglossarystyle{unifiedstyle}{%
  \setglossarystyle{list}%
  \renewcommand{\glossentry}[2]{%
    \edef\glscursec{\glsentryuseriii{##1}}%
    \ifx\glscursec\glslastsec\else
      \item[]\textbf{\S\,\glscursec\enspace\glsentryuseriv{##1}}%
      \global\let\glslastsec\glscursec
    \fi
    \item[\glsentryuserii{##1} \glsentryuseri{##1}:] ``\glossentryname{##1}'' \glsentrydesc{##1}%
  }%
}

\makeatletter
\def\@firstwordof#1 #2\relax{#1}

\newcommand{\beXa}[1][]{%
  \def\@temparg{#1}%
  \ifx\@temparg\@empty
    \begin{example}%
  \else
    \let\@savedSecTitle\@currentsectionname%
    \refstepcounter{excounter}%
    \begin{example}[#1]%
    \label{e:\theexcounter}%
    \expanded{%
      \noexpand\newglossaryentry{e\theexcounter}{%
        name={\unexpanded{#1}},%
        sort={\theexcounter},%
        description={p.~\noexpand\pageref{e:\theexcounter}},%
        user1={\theexcounter},%
        user2={Example},%
        user3={\thesection},%
        user4={\expandafter\unexpanded\expandafter{\@savedSecTitle}}%
      }%
    }%
  \fi
}
\newcommand{\eeXa}{\end{example}}

\newcommand{\beD}[1][]{%
  \def\@temparg{#1}%
  \ifx\@temparg\@empty
    \begin{definition}%
  \else
    \let\@savedSecTitle\@currentsectionname%
    \refstepcounter{defcounter}%
    \begin{definition}[#1]%
    \label{d:\thedefcounter}%
    \expanded{%
      \noexpand\newglossaryentry{d\thedefcounter}{%
        name={\unexpanded{#1}},%
        sort={\thedefcounter},%
        description={p.~\noexpand\pageref{d:\thedefcounter}},%
        user1={\thedefcounter},%
        user2={Definition},%
        user3={\thesection},%
        user4={\expandafter\unexpanded\expandafter{\@savedSecTitle}}%
      }%
    }%
  \fi
}
\newcommand{\eeD}{\end{definition}}

\NewDocumentCommand{\beT}{O{} O{} o}{%
  \def\@temparg{#1}%
  \ifx\@temparg\@empty
    \begin{theorem}%
  \else
    \let\@savedSecTitle\@currentsectionname%
    \refstepcounter{thmcounter}%
    \IfValueTF{#3}{%
      \edef\@thmlblkey{#3}%
    }{%
      \edef\@thmlblkey{\@firstwordof#1 \relax}%
    }%
    \def\@tempshort{#2}%
    \ifx\@tempshort\@empty
      \expandafter\xdef\csname @tshort@\@thmlblkey\endcsname{\@thmlblkey}%
    \else
      \expandafter\def\csname @tshort@\@thmlblkey\endcsname{#2}%
    \fi
    \begin{theorem}[#1]%
    \label{t:\@thmlblkey}%
    \expanded{%
      \noexpand\newglossaryentry{t\thethmcounter}{%
        name={\unexpanded{#1}},%
        sort={\thethmcounter},%
        description={p.~\noexpand\pageref{t:\@thmlblkey}},%
        user1={\thethmcounter},%
        user2={Theorem},%
        user3={\thesection},%
        user4={\expandafter\unexpanded\expandafter{\@savedSecTitle}}%
      }%
    }%
  \fi
}
\newcommand{\eeT}{\end{theorem}}
\newcommand{\beO}[1][]{%
  \def\@temparg{#1}%
  \ifx\@temparg\@empty
    \begin{open}%
  \else
    \let\@savedSecTitle\@currentsectionname%
    \refstepcounter{opcounter}%
    \begin{open}[#1]%
    \label{op:\theopcounter}%
    \expanded{%
      \noexpand\newglossaryentry{o\theopcounter}{%
        name={\unexpanded{#1}},%
        sort={\theopcounter},%
        description={p.~\noexpand\pageref{op:\theopcounter}},%
        user1={\theopcounter},%
        user2={Problem},%
        user3={\thesection},%
        user4={\expandafter\unexpanded\expandafter{\@savedSecTitle}}%
      }%
    }%
  \fi
}
\newcommand{\eeO}{\end{open}}
\newcommand{\beAs}[1][]{%
  \def\@temparg{#1}%
  \ifx\@temparg\@empty
    \begin{assumption}%
  \else
    \let\@savedSecTitle\@currentsectionname%
    \refstepcounter{ascounter}%
    \begin{assumption}[#1]%
    \label{as:\theascounter}%
    \expanded{%
      \noexpand\newglossaryentry{as\theascounter}{%
        name={\unexpanded{#1}},%
        sort={\theascounter},%
        description={p.~\noexpand\pageref{as:\theascounter}},%
        user1={\theascounter},%
        user2={Assumption},%
        user3={\thesection},%
        user4={\expandafter\unexpanded\expandafter{\@savedSecTitle}}%
      }%
    }%
  \fi
}
\newcommand{\eeAs}{\end{assumption}}
\newcommand{\beCj}[1][]{%
  \def\@temparg{#1}%
  \ifx\@temparg\@empty
    \begin{conjecture}%
  \else
    \let\@savedSecTitle\@currentsectionname%
    \refstepcounter{cjcounter}%
    \begin{conjecture}[#1]%
    \label{cj:\thecjcounter}%
    \expanded{%
      \noexpand\newglossaryentry{cj\thecjcounter}{%
        name={\unexpanded{#1}},%
        sort={\thecjcounter},%
        description={p.~\noexpand\pageref{cj:\thecjcounter}},%
        user1={\thecjcounter},%
        user2={Conjecture},%
        user3={\thesection},%
        user4={\expandafter\unexpanded\expandafter{\@savedSecTitle}}%
      }%
    }%
  \fi
}
\newcommand{\eeCj}{\end{conjecture}}
\newcommand{\beL}[1][]{%
  \def\@temparg{#1}%
  \ifx\@temparg\@empty
    \begin{lemma}%
  \else
    \let\@savedSecTitle\@currentsectionname%
    \refstepcounter{lemcounter}%
    \begin{lemma}[#1]%
    \label{l:\thelemcounter}%
    \expanded{%
      \noexpand\newglossaryentry{l\thelemcounter}{%
        name={\unexpanded{#1}},%
        sort={\thelemcounter},%
        description={p.~\noexpand\pageref{l:\thelemcounter}},%
        user1={\thelemcounter},%
        user2={Lemma},%
        user3={\thesection},%
        user4={\expandafter\unexpanded\expandafter{\@savedSecTitle}}%
      }%
    }%
  \fi
}
\newcommand{\eeL}{\end{lemma}}
\newcommand{\beR}[1][]{%
  \def\@temparg{#1}%
  \ifx\@temparg\@empty
    \begin{remark}%
  \else
    \let\@savedSecTitle\@currentsectionname%
    \refstepcounter{remcounter}%
    \begin{remark}[#1]%
    \label{rem:\theremcounter}%
    \expanded{%
      \noexpand\newglossaryentry{rem\theremcounter}{%
        name={\unexpanded{#1}},%
        sort={\theremcounter},%
        description={p.~\noexpand\pageref{rem:\theremcounter}},%
        user1={\theremcounter},%
        user2={Remark},%
        user3={\thesection},%
        user4={\expandafter\unexpanded\expandafter{\@savedSecTitle}}%
      }%
    }%
  \fi
}
\newcommand{\eeR}{\end{remark}}

\newcommand{\@currentsectionname}{}
\AtBeginDocument{%
  \let\@epidcrn@orig@refstepcounter\refstepcounter
  \renewcommand{\refstepcounter}[1]{%
    \@epidcrn@orig@refstepcounter{#1}%
    \ifstrequal{#1}{section}{\global\let\@currentsectionname\@currentlabelname}{}%
  }%
}

\newcommand{\beTa}[1][]{%
  \def\@temparg{#1}%
  \ifx\@temparg\@empty
    \begin{table}[htbp]%
  \else
    \let\@savedSecTitle\@currentsectionname%
    \refstepcounter{tabcounter}%
    \begin{table}[htbp]%
    \label{A:\thetabcounter}%
    \expanded{%
      \noexpand\newglossaryentry{A\thetabcounter}{%
        name={\unexpanded{#1}},%
        sort={\thetabcounter},%
        description={p.~\noexpand\pageref{A:\thetabcounter}},%
        user1={\thetabcounter},%
        user2={Table},%
        user3={\thesection},%
        user4={\expandafter\unexpanded\expandafter{\@savedSecTitle}}%
      }%
    }%
  \fi
}
\newcommand{\eeTa}{\end{table}}

\newcommand{\beF}[1][]{%
  \def\@temparg{#1}%
  \ifx\@temparg\@empty
    \begin{figure}[H]%
  \else
    \begin{figure}[H]%
    \let\@savedSecTitle\@currentsectionname%
    \refstepcounter{figcounter}%
    \label{fig:\thefigcounter}%
    \expanded{%
      \noexpand\newglossaryentry{Fig\thefigcounter}{%
        name={\unexpanded{#1}},%
        sort={\thefigcounter},%
        description={p.~\noexpand\pageref{fig:\thefigcounter}},%
        user1={\thefigcounter},%
        user2={Figure},%
        user3={\thesection},%
        user4={\expandafter\unexpanded\expandafter{\@savedSecTitle}}%
      }%
    }%
  \fi
}
\newcommand{\eeF}{\end{figure}}
\makeatother
\newcommand{\printboth}{%
  \clearpage
  \glsaddall  
  \printglossary[style=unifiedstyle, title={0 \hspace{.22em} Definitions, Theorems, Examples and Open Problems}]%
}
\newcommand{\tref}[1]{%
  Thm~\ref{t:#1}(\ifcsname @tshort@#1\endcsname\csname @tshort@#1\endcsname\else#1\fi)%
}

\makeatletter
\@ifundefined{lemma}      {\newtheorem{lemma}      {Lemma}}      {}
\@ifundefined{proposition}{\newtheorem{proposition}{Proposition}} {}
\@ifundefined{corollary}  {\newtheorem{corollary}  {Corollary}}  {}
\@ifundefined{conjecture} {\newtheorem{conjecture} {Conjecture}}  {}
\@ifundefined{remark}     {\newtheorem{remark}     {Remark}}     {}
\@ifundefined{question}   {\newtheorem{question}   {Question}}   {}
\@ifundefined{assumption} {\newtheorem{assumption} {Assumption}} {}
\@ifundefined{hypothesis} {} {}
\@ifundefined{fact}       {}       {}
\makeatother

\def\beP{\begin{proposition}}\def\eeP{\end{proposition}}
\def\beC{\begin{corollary}}\def\eeC{\end{corollary}}
\def\beQ{\begin{question}}\def\eeQ{\end{question}}
\def\beA{\begin{assumption}}\def\eeA{\end{assumption}}

\definecolor{funccolor}{RGB}{25,25,112}
\definecolor{desccolor}{RGB}{64,64,64}

\def\eq{=}

\def\bep{\begin{pmatrix}}\def\eep{\end{pmatrix}}
\def\bev{\begin{vmatrix}}\def\eev{\end{vmatrix}}
\def\bea{\begin{eqnarray*}}\def\eea{\end{eqnarray*}}
\def\bc{\begin{cases}}\def\ec{\end{cases}}
\def\BEN{\begin{enumerate}}\def\EEN{\end{enumerate}}
\def\BI{\begin{itemize}}\def\EI{\end{itemize}}

\newcommand{\be}[1]{\begin{equation}\label{#1}}
\newcommand{\ee}{\end{equation}}
\newcommand{\beq}{\begin{eqnarray}}
\def\eeq{\end{eqnarray}}
\def\eqr{\eqref}\def\lbl{\label}
\def\Lra{\Longrightarrow}

\newcommand{\R}{\mathbb{R}}  

\def\T{\widetilde}

\newcommand{\bff}[1]{{\mbox{\boldmath$#1$}}}

\newcommand\y{\boldsymbol{y}}

\def\v1{\vec {\bff 1}}

\def\mL{\mathcal{L}}
\def\mR{{\mathcal R}}

\newcommand\op{\operatorname}
\newcommand{\supp}{\op{supp}}

\long\def\symbolfootnote[#1]#2{%
\begingroup
\def\thefootnote{\fnsymbol{footnote}}\footnote[#1]{#2}%
\endgroup}

\def\and{antisymmetric}

\def\ie{i.e. }
\def\im{\item}

\def\ssec{\subsection}

\def\var{variable}

\newcommand\CRN{chemical reaction network}

\def\brn{basic reproduction number}

\def\DFE{disease free equilibrium}

\def\ME{mathematical epidemiology}

\def\NGM{next generation matrix}

\def\fp{fixed point}

\def\Lf{Lyapunov function}

\def\regS{{\bf regular splitting}}

\def\NGMs{next generation matrices}

\def\PV{Perron-Volterra candidate Lyapunov function}

%% file: intR1.tex
\section{Introduction}
\subsection{Short review of some Volterra  type Lyapunov function results in ME and CRN}\lbl{s:Lya}

Since Volterra, Lyapunov  functions for positive ODEs $z'=f(z)$ have been often constructed
using the Bregman divergence:
\[
  \T G_i(z_i):=z_i-z_i^*-z_i^*\,\ln\frac{z_i}{z_i^*}=z_i^* G(u_i)\ge 0,
  \qquad u_i:=\frac{z_i}{z_i^*}, G(u)=u-1-ln(u)
\]
 associated with $x\mapsto x\ln x$.
 ``Candidate Lyapunov functions" include often terms of the form:
\be{Vol}
  \mathcal L(z):=\sum_{i=1}^n l_i \T G_i(z_i),\qquad l_i>0.
\ee

Note that $\T G_i'=(1-u_i^{-1})z_i'$ yields the explicit formula:
\begin{equation}\label{eq:Lprime_raw}
  \mathcal L'
  = \sum_{i=1}^n l_i\!\Bigl(1-\frac{1}{u_i}\Bigr)f_i(z).
\end{equation}

\paragraph{Unconditional GAS under complex balance and strong endotacticity in \CRN.}
Classical CRNT \cite{HornJackson} proves \emph{unconditional} global asymptotic stability (GAS) of  complex-balanced equilibria within their stoichiometric compatibility class, {\bf for arbitrary positive rate constants}, by means of the Volterra--Horn--Jackson  Lyapunov function
\[
  \mL(z)=\sum_i \bigl(z_i\ln(z_i/z_i^*)-z_i+z_i^*\bigr)
\]
(see also \cite{AndersonGAS} for an extension of weakly reversible, single
 linkage mass-action ODEs,
 and see \cite{GopMilShiu}
for an extension to strong endotactic networks).
 Note here the Lyapunov argument is purely structural: once complex balancedness or strong endotacticity  hold, the construction of the \Lf\ does not require any computation, beyond that of a positive steady state.

 \paragraph{Conditional GAS in \ME\ models.}
 Fall, Iggidr, Sallet, and Bonzi \cite{Fall,IggidrCEP,Bonzi}--
 see also Earn and McCluskey \cite{EarnMc}
and~\cite{AAHK}---have studied a class of bilinear epidemic models
with rank-one next-generation matrix (NGM), in which the right Perron
eigenvector yields an explicit endemic equilibrium~(EE), while the left
Perron eigenvector provides Lyapunov functions at both the disease-free
equilibrium~(DFE) and the~EE.
Furthermore, they established strong
threshold theorems in the sense of Shuai and van den Driessche (2013):
when the DFE loses stability, a unique EE emerges and becomes globally
asymptotically stable, with the transition certified by an explicit global
Lyapunov function. For the non-rank one case,  Shuai and van den Driessche \cite{Shuai13} presented a  method for   proving the global stability of the disease-free equilibrium (DFE), using a Lyapunov function which   combined  a Volterra function in the variable(s) which are unconditionally positive, and a linear function in the ``infection/siphon/possibly zero" \var s, with weights given by  the left Perron eigenvector of  the \NGM\ (NGM). A  second result of \cite{Shuai13} pertained to the  more difficult problem of establishing \emph{conditional} GAS for a selected endemic equilibrium (EE). This second result
is however not fully constructive, depending on the ability of the user to identify a graph structure, and to establish certain inequalities.

A striking further result was obtained by~\cite{Rahman}, who proved a strong
threshold theorem for a two-strain model with rational saturating rates.  In
this setting, stability is relayed across equilibria with decreasing infection
levels, inducing a partition of parameter space into regions where exactly one
equilibrium is certified to be globally stable by an explicit Lyapunov
function.  We refer to this as the competitive exclusion partition
property~(CEPP).

Our work is an extension of the results of \cite{Fall,IggidrCEP,Bonzi,EarnMc} to models with several boundary equilibria like in \cite{Rahman}, and also by considering  increasing concave incidence functions --see Theorem~\ref{t:GAS}. We also show,  by providing explicit global Lyapunov functions
at all equilibria--see
Theorem~\ref{t:blo}, that competitive exclusion partitions of the parameter space hold for models with $k$ strains, as long as all
strains are all rank-one and non-interacting.
The Lyapunov functions are sums of a
Volterra entropy on resident variables and a linear functional on invaders,
weighted by left Perron eigenvectors of transversal Jacobians.

Besides unifying and extending explicit Lyapunov stability results
of~\cite{Fall,IggidrCEP,Bonzi,EarnMc} and~\cite{Rahman}, we also reveal their connection to CRNT, and implement
the results algorithmically in the Mathematica package \texttt{EpidCRN}, which
views multi-strain epidemic models through the lens of chemical reaction
network (CRN) theory.

\beR[Combining the Rahman--Zou complete partition GAS result with the CRN structural approach]

The remarkable competitive exclusion partition obtained by Rahman and Zou
(and which persists even in the presence of delays) raises the natural
question of identifying the largest class of epidemic models for which
such a complete partition remains valid.

The CRN structural viewpoint of modifying the model suggests that the key object is not the particular form
of the incidence functions, but the rank-one structure of the Metzler transversal
Jacobian $J^\perp_0$ at the strain siphon faces.  

Motivated by a  structure search, we obtained below several extensions of the
Rahman--Zou theorem, including models with increasing concave incidence,
and some $n$-strain  rank-one models.

One structural assumption that remains unstudied concerns the
Frobenius graph associated with $J^\perp_0$.  All results obtained here
correspond to the case where this graph consists of $n$ isolated
vertices, equivalently, where the Frobenius decomposition is block
diagonal and the irreducible components are dynamically independent.
Whether a complete exclusion partition continues to hold for some
Frobenius forests with nontrivial couplings between distinct components
remains open -- see
Problem~\ref{o:GASrkone-general}.

\eeR

%% file: mt.tex
\subsection{The m-strain model with  rank one invading blocks}
One of the most convincing results for the utility of CRN-ME tools is a
result for the  model with $m$ minimal disjoint siphons below.

Consider $m$ competing irreducible
rank-one blocks, with no direct block-to-block coupling in the
transversal equations,
\be{eq:rankone-m}
\begin{cases}
s'
=
\Lambda-\mu s-s\displaystyle\sum_{k=1}^m\phi_k,
\\[1mm]
z_k'
=
\left(sw_k\pi_k^\top-V_k\right)z_k,
\qquad
k=1,\ldots,m,
\\[1mm]
\phi_k
=
\pi_k^\top z_k.
\end{cases}
\ee
Here $z_k\in\R_{\ge0}^{n_k}$, $w_k\gg0$, $\pi_k\gg0$,
${\bf 1}^\top\pi_k=1$, and    $V_k$ are  nonsingular
$M$-matrices, with  $V_k^{-1}\gg0$ (internal transitions/removal).

Note the ``structure aware" writing of the ODE, which partitions the variables into  $m$ non-interacting disjoint minimal siphon/maximal invariant faces/``generalzed Lotka-Volterra"/invading  groups, and a non-siphon/non-DFE variable $s$, and the rank-one regular splitting of the transversal Jacobian/invading matrices for each block:
\[
J_k^\perp(s)=F_k(s)-V_k,
\qquad
F_k(s)=sw_k\pi_k^\top\ \gg0(\text{rank one, new infections}).
\]

We assume that 
\be{wr}
{\bf 1}^\top V_k\geq\mu{\bf 1}^\top,
\qquad
k=1,\ldots,m
\ee
(in epidemiologic terms, the decay of the invading strain compartments surpasses
that of the input, ensuring the invariance of the set below).
Then
\[
\Omega
=
\left\{
s\geq0,\ z_k\geq0:
s+\sum_{k=1}^m{\bf 1}^\top z_k\leq s_0
=\frac{\Lambda}{\mu}
\right\}
\]
is forward invariant.

%% file: Fr.tex
\paragraph{Frobenius normal form.}
We pause to recall a useful result on the spectral structure of Metzler
matrices~\cite{HornJohnson}:

\beP[Frobenius normal form of Metzler matrices]\lbl{p:FNF}
Let $M\in\R^{n\times n}$ be Metzler (i.e.\ $M_{ij}\ge 0$ for $i\ne j$).
Then there exists a permutation matrix $P$ such that
\[
\T M=P^{-1}MP=
\begin{pmatrix}
\T B_1 & * & \cdots & *\\
0   & \T B_2 & \cdots & *\\
\vdots & & \ddots & *\\
0   & 0   & \cdots & \T B_r
\end{pmatrix},
\]
where each diagonal block $\T B_j$ is an {\bf irreducible} Metzler matrix.
Moreover:
\begin{enumerate}
\item $\lambda_{\max}(M)=\max_{1\le j\le r}\lambda_{\max}(\T B_j)$, where
  $\lambda_{\max}$ denotes spectral abscissa.

\item For each $j$, the Perron--Frobenius theorem yields $u_j\gg0$ with
  $\T B_j u_j=\lambda_{\max}(\T B_j)u_j$.

\item If the form is \emph{block diagonal} (all $*=0$), then extending $u_j$ by zeros and
conjugating by $P$ gives a nonnegative right eigenvector of $M$ at $\lambda_{\max}(\T B_j)$; in
particular the block $\T B_{j^*}$ attaining $\lambda_{\max}(M)$ yields the Perron eigenvector of
$M$. In the general upper-triangular case the Perron eigenvector of $M$ need not be supported
on a single Frobenius block, since a zero-extension of $u_j$ satisfies
$(P^{-1}MP)$-eigenvector equation only when the off-diagonal blocks above $j$ annihilate
$u_j$.
\end{enumerate}

\eeP

%% file: ct.tex
\ssec{Contents}    The paper contains:
\BEN
\im A  section~\ref{s:def}, where some CRN and ME definitions are provided.
\im A  section~\ref{s:T1}, where a multi-strain GAS theorem is  stated.

\im A competitive exclusion  GAS-CEP partition -- see theorem ~\ref{t:blo}-- for either models with a finite collection of irreducible  non-interacting
rank-one strains and diagonal matrices $V_k$, or for two strain models with arbitrary $V_k$, in section~\ref{s:blo}.
\im
 A  proof of a competitive exclusion partition   for
two-strain concave incidence models (section~\ref{s:Rah}),  generalizing \cite{Rahman}, and   providing explicit
global Lyapunov functions at all four equilibria via a single Perron--Volterra ansatz.


\im Section~\ref{s:alg} illustrates how Epid-CRN helps with solving one example LyapExa.wl.

\im Section~\ref{s:con} concludes and formulates open problems.

\EEN

%% file: wh.tex
\section{Why combine CRN and ME methods?}\lbl{s:def}
The first characteristic of most models in mathematical epidemiology,
ecology, and immuno-virology is that they are positive ODEs (or stochastic
models whose mean field is a positive ODE).  This means that important
phenomena occur when trajectories reach the boundary of the positive orthant,
and that the notion of invariant boundary faces is crucial.  While boundary
faces have been studied in classical ODE theory, the systematic exploitation
of simplifications arising for positive ODEs was initiated by CRN works; see
notably~\cite{hun,Erdi}, who showed that all positive polynomial models admit
(non-unique) mass-action representations.

Beyond the polynomial case, but assuming a stoichiometric representation,
crucial results relating boundary $\omega$-limit points to invariant boundary
faces and to the Petri-net concept of siphons were obtained in~\cite{AdLS}.
A recent proof of the folklore result that the Jacobian has a triangular block
form on siphon faces was given in~\cite{AAH26,AH26}, and \cite{AAHK} noted
that transversal Jacobians on siphons have the Metzler/cooperativity property.
This property further ensures the existence of left and right Perron
eigenvectors which, in all cases we have studied, provide Lyapunov functions
and escape directions from siphon faces respectively.

Our framework is encapsulated by the following definitions.

\beD[Positive / non-negative ODE]
An ODE is called \emph{positive}~\cite{rantzer2015scalable} or
\emph{non-negative}~\cite{haddad2010} if the non-negative orthant
\[
\mathbb{R}^n_{\ge 0} := \{ x \in \mathbb{R}^n : x_i \ge 0,\; i = 1,\dots,n \}
\]
is forward invariant under the flow.
\eeD

\beD[Stoichiometric and chemical representation]\label{d:crn}
A stoichiometric representation of $f$ is a pair $(\Gamma, r)$ such that
\[
f(x) = \Gamma r(x), \qquad r(x) \ge 0,
\]
where $\Gamma \in \mathbb{R}^{n \times n_R}$ is constant and
$r : \mathbb{R}^n_{\ge0} \to \mathbb{R}^{n_R}_{\ge0}$ is locally Lipschitz.

An ODE is called chemical if
\[
\Gamma_{i\rho} < 0 \implies r_\rho(x) = x_i \tilde{r}_\rho(x), \qquad
\tilde{r}_\rho(x) \ge 0,
\]
and, moreover, $r_\rho$ depends only on variables in
$\supp(r_\rho) \subseteq \Gamma^-_i$ whenever $\Gamma_{i\rho} < 0$.
Furthermore, we assume that each rate $r_\rho(x)$ is monotone non-decreasing
with respect to each variable $x_k$ for $k \in \supp(r_\rho)$.
\eeD

\paragraph{Siphons.}
Angeli, de Leenheer and Sontag~\cite{AdLS}, Anderson~\cite{AndGAS} and Shiu
and Sturmfels~\cite[Prop.~2.1]{ShiuStu} proved that a boundary face $F_\sigma$
associated to a nonempty set $\sigma$ of variables being $0$ is
forward-invariant for a chemical ODE iff $\sigma$ is a siphon in the sense of
the combinatorial Definition~\ref{d:sip}:

\beD[Siphon, minimal siphon, total siphon/DFE, critical siphon]\lbl{d:sip}
\leavevmode
\begin{itemize}
\item A \textbf{siphon set} $\sigma \subseteq \mathcal{S}$ is a nonempty
  subset of species such that whenever a species in $\sigma$ appears in a
  product complex in the RHS of a reaction, at least one species in $\sigma$
  must appear also in the LHS of the reaction (the reactant complex).

\im A siphon is \textbf{minimal} when it contains no other siphon.

\im The union of all minimal siphons will be called the \textbf{total
  siphon/DFE} and denoted~$\sigma_0$.

\im A siphon $\sigma$ in a CRN with stoichiometric matrix $\Gamma$ is
  \textbf{critical} if it contains no support of a positive conservation
  relation, i.e., if there exists no nonzero vector $c \ge 0$ with
  $c \Gamma=0$ and $\operatorname{supp}\, c \subset \sigma$.
\end{itemize}
\eeD

\beR[CRNT suggests definitions for important ME concepts]
The total siphon/DFE set $\sigma_0$ is a fundamental object of study in
mathematical epidemiology (ME), but a rigorous definition for it is available only
via CRNT.  Similarly, CRNT suggests definitions for multi-strain and ME models--see below.
\eeR

We have now the building bricks to define ME and multi-strain systems.
\beD[ME models, $k$--siphon  models, $k$--strain simple models]
\lbl{d:ME}
Consider a chemical ODE with at least one critical minimal siphon, and a unique DFE, and let
$J^\perp_{0}$ denote the diagonal Jacobian block in the triangularization of
the Jacobian corresponding to the DFE siphon, which exists and is Metzler
by~\cite{AAH26,AH26,AAHK}.

 An \textbf{epidemic model} is a pair formed from  a chemical ODE and
 a regular splitting~(RS) $J^\perp_{0}=F_0-V_0$ of the transversal Jacobian $J^\perp_{0}$ in
the sense of~\cite{Varga,Van} (with $F_{0}\geq 0$ and $V_{0}$ a non-singular
M-matrix, hence $V_0^{-1}\geq 0$). Note that since $J^\perp_{0}$ is Metzler \cite{AAH26}, a trivial \regS\ always exists if  all the diagonal elements of $J^\perp_{0}$ have strictly negative parts.

A $k$--siphon  model is an ME model with precisely $k$ minimal siphons.

A $k$--strain simple model is a chemical ODE
with precisely $k$ disjoint minimal siphons $\sigma_1,\dots,\sigma_k$.
\eeD

\beT[Generalized Lotka--Volterra form of $k$--strain simple models]
A $k$--strain simple model with variables
\[
(s,z^1,\dots,z^k)
\in
\mathbb R_{\ge0}^m
\times
\prod_{j=1}^k\mathbb R_{\ge0}^{d_j},
\] where each block $z^j$ represents the variables/infected class associated with
the strain/siphon $\sigma_j$, and
$s=(s_1,\dots,s_m)$ represents the susceptible or uninfected
classes,
may always be written in generalized
Lotka--Volterra form, with the minimal-siphon equations factored, and ``regularly split", as follows:
\be{MS}
\dot s
=
f_s(s,z^1,\dots,z^k),
\qquad
\dot z^j
=
\bigl(F_j-V_j\bigr)(s,z^1,\dots,z^k)z^j,
\qquad
j=1,\dots,k.
\ee

Here, for each $j$, letting $\y_j$ denote the variables which remain free on the
face $z^j=0$, and putting
\(
F_j(\y_j)
=
F_j(s,z^1,\dots,z^k)\big|_{z^j=0},
V_j(\y_j)
=
V_j(s,z^1,\dots,z^k)\big|_{z^j=0},
\)
we may  chose $F_j, V_j$ so that
\[
J^\perp_{\sigma_j}(\y_j)
=
F_j(\y_j)-V_j(\y_j)
\]
is a regular splitting (i.e. $F_j(\y_j)\ge 0, V_j(\y_j)^{-1}\ge 0$).

\eeT

\begin{proof}
Since $\sigma_j$ is a siphon, the face $z^j=0$ is forward invariant;
hence the $z^j$-component $g_j$ of the vector field satisfies
$g_j(s,z^1,\dots,0,\dots,z^k)=0$.  By smoothness,
\[
g_j(s,z^1,\dots,z^k)
=
A_j(s,z^1,\dots,z^k)\,z^j,
\qquad
A_j=
\int_0^1
D_{z^j}g_j(s,z^1,\dots,tz^j,\dots,z^k)\,dt,
\]
which yields =ref{MS} with $A_j$ in place of $F_j-V_j$. Setting $t=0$
in the integrand shows $A_j(\y_j,0)= D_{z^j}g_j(s,z^1,\dots,0,\dots,z^k)
= J^\perp_{\sigma_j}(\y_j)$, so $A_j$ restricted to the face $z^j=0$
is exactly the transversal Jacobian.

Since $J^\perp_{\sigma_j}(\y_j)$ is Metzler by \cite{AAH26}, it admits
some regular splitting $J^\perp_{\sigma_j}(\y_j)= F_j^0(\y_j)-V_j^0(\y_j)$
(the trivial diagonal one always works). It remains to extend $F_j^0,V_j^0$
off the face so that $A_j= F_j-V_j$ holds identically and the
extensions restrict to $F_j^0,V_j^0$ on $z^j=0$. Define
\[
F_j(s,z^1,\dots,z^k)= F_j^0(\y_j)
\]
(i.e.\ $F_j$ is simply $F_j^0$, frozen at its face value and read off the
$z^j$-independent variables $\y_j$; in particular $F_j\ge0$ everywhere,
since $F_j^0\ge0$), and
\[
V_j= F_j-A_j.
\]
Then $A_j= F_j-V_j$ holds everywhere by construction, and on the face,
$F_j|_{z^j=0}= F_j^0$ while
$V_j|_{z^j=0}= F_j^0-A_j(\y_j,0)= F_j^0-\bigl(F_j^0-V_j^0\bigr)= V_j^0$,
so this simultaneous choice recovers exactly the originally chosen
regular splitting on the face, as claimed. This completes the construction: $F_j,V_j$ have been
exhibited explicitly, for every $j$, from the data of the chemical ODE and the chosen face
splitting $F_j^0-V_j^0$ alone.
\end{proof}

\beR[k-strain simple models generalize Lotka-Volterra models]
The generalization arises from two directions.  First,  not all equations are required to have multiplicative structure. Secondly, the scalar multiplications
in Lotka-Volterra (corresponding to scalar siphons) are replaced by matrix multiplications in =r{MS}. \eeR

\beD[$R$--invasion functions and rank-one $k$--strain models]\label{d:RepFun}
For a simple k-strain model, for each strain $j$ and regular splitting, define \NGMs\ $K_j(\y_j)=(F_j V_j^{-1})(\y_j)$ and invasion functions
\[
 R_j(\y_j):=\rho\bigl(K_j(\y_j)\bigr),
\]
 $R_j(\y_j)$ will be called {\bf R-invasion functions}
(associated with the  block $j$).

A \textbf{simple k-strain model is rank-one} if every Jacobian transversal block $J^\perp_{\sigma}$
 admits a  regular splitting $J^\perp_{\sigma}=F_\sigma-V_\sigma$, such that the  ``new infections" matrix $F_\sigma$ has rank one. In this case, the  R-invasion functions and the Perron data admit explicit rank-one formulas.
\eeD

\beR[Previous appearances of the R-invasion functions]
Interestingly, invasion functions appear in a different context in the work of Magal, Webb, and Wu (2019) \cite{MagWW}, who develop  threshold-based methods for nonlinear epidemic models in a spatial reaction--diffusion setting.
\eeR

\beD[reproduction  numbers and invasibility numbers] \label{d:reF}
\leavevmode

\begin{enumerate}
\im  The evaluations
\be{Ri} R_i=R_i(E_0), i=1,...,k,\ee
where $E_0$ is the DFE, and $R_i(\y_i)$ are the invasion functions, will be called   {\bf \brn s} (of strain $i$).

\im   In the case  when each minimal siphon  has a unique resident \fp\
$E_j, j=1,...,k,$
\be{invN}R_i^{\T j}:=R_i(E_j)\ee will be called
the {\bf invasibility number of invading block $i$ on resident block $j$}.

\EEN
\eeD

A central motivation of our paper is the observation that the explicit
results of~\cite{Fall,IggidrCEP,Bonzi,EarnMc,Shuai13} and others share a
common structure: they are rank-one models, and exploit the computation of
left and right Perron eigenvectors of rank-one NGMs.

\beR[The roles of the Perron eigenvectors of transversal Jacobians on siphon faces]
Once Metzler transversal Jacobians on siphon faces have been identified as
the basic objects governing invasion, it is natural to study the roles of
their Perron eigenvectors.

The right Perron vector determines the invasion direction. In the rank-one
models considered in this paper and in  rank one works like for example \cite{Fall,IggidrCEP,Bonzi,EarnMc,AAHK}, it also yields the coordinates of the
boundary  equilibria.

The left Perron vector plays a dual role. It determines both the canonical
normalized vector
\[
\T \pi_k^\top
=
\frac{\pi_k^\top V_k^{-1}}{\mR_k},
\]
used in the linear invader part of the Perron--Volterra Lyapunov function,
and the associated resident entropy weights
\[
\T w_{k,j}=\T \pi_{k,j}w_{k,j},
\]
normalized so
\(
\sum_j\T w_{k,j}=1.
\)
\eeR

\beR[The choice of Frobenius decomposition and splitting fixes the  Perron data]
\lbl{r:canonical}

Throughout the paper, we consider   epidemic models (Definition~\ref{d:ME}),
with a fixed regular splitting
\[
J_{\sigma_0}^\perp(E_0)= F_0-V_0,
\]
and the Frobenius decomposition and   splittings are fixed once and for
all, fixing the Perron  data.
\eeR

%% file: GASt.tex
\section{A GAS theorem~\ref{t:blo}  for some rank-one simple multi-strain models}\lbl{s:T1}
One of the most convincing results pointing to the interest of
ME--CRN unification is the GAS competitive exclusion partition
theorem for rank-one multi-strain models.
\subsection{Two GAS results for  rank-one simple multi-strain models}
 The  theorem below lists the equilibrium and invasion data
for the  generic equilibrium types, and states extra conditions under which a GAS CEP holds, in two cases; the degenerate coexistence
case is treated separately in Remark~\ref{r:generic-rankone} below.

\beT[Global rank-one blocks competition theorem]
\lbl{t:blo}
For the rank-one m-strain model \eqref{eq:rankone-m}, for each block, let
\[
\mR_k
=
\pi_k^\top V_k^{-1}w_k,
\qquad
R_k(s)
=
s\mR_k,
\qquad
R_k
=
s_0\mR_k,
\]
denote the per-susceptible reproduction coefficients, the  reproduction
functions and
the  reproduction numbers.

This model may have only generic equilibria supported on one or none rank-one blocks. Their equilibrium formulas and invasion numbers are
determined as follows.
\BEN \im
\BEN

\im
The disease-free equilibrium is
\[
E_0
\eq
(s_0,0,\ldots,0).
\]
The transversal eigenvalue of block $k$ at $E_0$ has the sign of
\[
R_k-1.
\]

\im
For a single block $p$, a boundary equilibrium supported only on block $p$
exists if and only if
\[
R_p>1.
\]
It is
\[
E_p
\eq
(s_p,0,\ldots,0,z_p^*,0,\ldots,0),
\]
where
\[
s_p
\eq
\frac1{\mR_p},
\qquad
z_p^*
\eq
\xi_pV_p^{-1}w_p,\qquad
\xi_p
\eq
\mu(s_0- s_p)
\eq
\Lambda(1-\frac{1}{R_p}).
\]
For $k\ne p$, the invasion number of block $k$ at $E_p$ is
\[
R_k(E_p)
\eq
s_p\mR_k
\eq
\frac{\mR_k}{\mR_p}
\eq
\frac{R_k}{R_p}.
\]
Hence block $k$ invades $E_p$ if and only if
\[
R_k>R_p,
\]
equivalently $\mR_k>\mR_p$.

\EEN

\im
Assume in addition that every resident block $p\in\{1,\ldots,m\}$ satisfies the
resident entropy closing condition of Definition~\ref{d:closing}; by
Lemma~\ref{l:two-mechanisms} this holds, block by block, whenever each
block $k\in\{1,\ldots,m\}$ individually satisfies either (D) $V_k$
diagonal, or (T) $n_k=2$ with $V_k$ an irreducible nonsingular $M$-matrix
(equivalently, satisfying the standing assumption $V_k^{-1}\gg0$ of
\eqref{eq:rankone-m})
-- the two cases may be mixed freely across blocks, e.g.\ some blocks
scalar/diagonal and others two-dimensional with an irreducible $M$-matrix
$V_k$. The
following global stability partition then holds generically (away from
the non-generic coexistence stratum of Remark~\ref{r:generic-rankone}).

\BEN

\im
If
\[
R_k\leq1,
\qquad
k\eq1,\ldots,m,
\]
then
\[
E_0
\eq
(s_0,0,\ldots,0)
\]
is globally asymptotically stable on $\Omega$.

\im
Suppose that, for some $p$,
\[
R_p>1,
\qquad
R_p>R_k
\quad
(k\ne p).
\]
Then the single-block equilibrium
\[
E_p
\eq
(s_p,0,\ldots,0,z_p^*,0,\ldots,0),
\]
where
\[
s_p
\eq
\frac1{\mR_p},
\qquad
\xi_p
\eq
\Lambda-\mu s_p
\eq
\Lambda-\frac{\mu}{\mR_p},
\qquad
z_p^*
\eq
\xi_pV_p^{-1}w_p,
\]
is globally asymptotically stable on the persistence class
\[
\Omega_p^+
\eq
\{(s,z_1,\ldots,z_m)\in\Omega:z_p\ne0\}.
\]

\im
On the non-generic tie surface of Remark~\ref{r:generic-rankone}, the
corresponding resident equilibrium is not isolated, and no point of
the continuum is GAS in the ordinary sense.
\EEN
\EEN
\eeT

\beR[Genericity, and the non-generic coexistence stratum]
\lbl{r:generic-rankone}

The equilibrium structure of \eqref{eq:rankone-m} is determined
entirely by the $m$ invasion functions $R_k(s)=s\mR_k$.  Each
positive block forces $s=1/\mR_k$, so a putative equilibrium with
support $I\subset\{1,\ldots,m\}$, $|I|\ge2$, requires
\[
\mR_k\eq\mR_I
\qquad\text{for every }k\in I
\]
for some common value $\mR_I$ -- that is, two or more of the $m$ real
numbers $\mR_1,\ldots,\mR_m$ must coincide exactly.  For generic
parameter values $\mR_1,\ldots,\mR_m$ are pairwise distinct, since
coincidence is a codimension-$(|I|-1)$ condition on the parameter
space.  Consequently, \emph{coexistence of two or more rank-one
blocks is non-generic}: for generic parameters the only equilibria of
\eqref{eq:rankone-m} are
\[
E_0,E_1,\ldots,E_m,
\]
as given by Theorem~\ref{t:blo}, and the GAS partition
theorem below (Theorem~\ref{t:blo}) settles the dynamics
completely in this generic case: either $E_0$ or exactly one $E_p$ is
GAS.

On the non-generic tie surface where $\mR_k=\mR_I$ for $k\in I$, a
short computation -- repeating the argument above for each $k\in I$
simultaneously -- shows that the equilibrium is not isolated: writing
$s=1/\mR_I$, $\xi_I=\Lambda-\mu/\mR_I$, the set of equilibria
supported on $I$ is the simplex
\[
E_I(\xi)
\eq
\left(
\frac1{\mR_I},
z_1(\xi),\ldots,z_m(\xi)
\right),
\qquad
z_k(\xi)
\eq
\begin{cases}
\xi_kV_k^{-1}w_k, & k\in I,\\
0, & k\notin I,
\end{cases}
\qquad
\xi_k>0,\ \sum_{k\in I}\xi_k\eq\xi_I,
\]
with invasion number $R_q(E_I)=\mR_q/\mR_I=R_q/R_I$ at any point of
the continuum, for $q\notin I$.  Since this stratum has codimension
$\ge1$ and supports a continuum (rather than an isolated point) of
equilibria, no point of it can be GAS in the ordinary sense; we do not
pursue its stability further here, and Theorem~\ref{t:blo}
below makes no GAS claim on this stratum, recording only that the
tied equilibrium is non-isolated.
\eeR

\begin{proof} Part 1:
at $E_0$, the transversal block for $z_k$ is
\[
s_0w_k\pi_k^\top-V_k.
\]
Since
\[
V_k^{-1}s_0w_k\pi_k^\top
\]
has spectral radius
\[
s_0\pi_k^\top V_k^{-1}w_k
\eq
s_0\mR_k
\eq
R_k,
\]
the sign of the Perron eigenvalue of
\[
s_0w_k\pi_k^\top-V_k
\]
is the sign of
\[
R_k-1.
\]

Now consider an equilibrium with block $p$ positive and all other blocks
zero. The equation
\[
0
\eq
(sw_p\pi_p^\top-V_p)z_p
\]
with $z_p>0$ implies, writing $\T i_p:=\pi_p^\top z_p$ for the (as yet unknown)
Perron aggregate infective level of block $p$,
\[
z_p
\eq
s\T i_pV_p^{-1}w_p.
\]
Applying $\pi_p^\top$ gives
\[
\T i_p
\eq
s\T i_p\pi_p^\top V_p^{-1}w_p
\eq
s\T i_p\mR_p.
\]
Since $\T i_p>0$, this gives
\[
s
\eq
\frac1{\mR_p}.
\]
The $s$-equation then gives
\[
0
\eq
\Lambda-\frac{\mu}{\mR_p}-\frac1{\mR_p}\T i_p,
\]
hence
\[
\T i_p
\eq
\mR_p\left(\Lambda-\frac{\mu}{\mR_p}\right)
\eq
\mR_p\xi_p.
\]
Thus
\[
z_p
\eq
s\T i_pV_p^{-1}w_p
\eq
\frac1{\mR_p}\mR_p\xi_pV_p^{-1}w_p
\eq
\xi_pV_p^{-1}w_p
\eq:
z_p^*,
\]
so that, equivalently,
\be{Vpzpstar}
V_pz_p^*\eq\xi_pw_p,
\ee
an identity used repeatedly below (Lemma~\ref{l:can-nor}, and the proof of
Theorem~\ref{t:n2-closing}). This equilibrium is positive exactly when $\xi_p>0$, equivalently
\[
\Lambda-\frac{\mu}{\mR_p}>0,
\]
equivalently
\[
s_0\mR_p>1,
\]
that is,
\[
R_p>1.
\]

At $E_p$, the missing block $k$ sees susceptible level $s_p$. Hence its
invasion number is
\[
R_k(E_p)
\eq
s_p\mR_k
\eq
\frac{\mR_k}{\mR_p}
\eq
\frac{R_k}{R_p}.
\]

The proof of the second (GAS-CEP) part of the theorem, based on establishing the decreasingness of a \PV,  is postponed  for section \ref{s:1}.
\end{proof}

%% file: not.tex
\subsection{Basic concepts, notations for defining  \PV s, for rank-one models}
For  rank one models, we will obtain below GAS partitions of the parameter space using Lyapunov functions \eqr{Hp} which combine Volterra entropy terms for the resident variables with Perron-weighted linear functionals for the invading variables, whose weights  obtained from the rank one decompositions of the  new infections matrices $F_\sigma$. Defining the family of \PV s for this model   uses several  concepts: 
\BEN \im      siphons  (a \CRN\ and Petri nets concept), which yield forward‑invariant coordinate faces for chemical ODEs;
\im  the disease‑free siphon, which is  the union of minimal
 siphons;
  \im the fact that Jacobians on siphon faces have a triangular block form, cf.
 \cite{AAH26}; \im the  fact  that a transversal Jacobian block on a siphon face is Metzler \cite{AAHK}, which puts under spotlight the roles of its Perron eigenvectors;
 \im the Perron-Frobenius decomposition of the transversal Jacobian on the DFE siphon face $J^\perp_{0}$, which defines the ``invading strains", with transversal blocks $J^\perp_{\sigma}$;
 \im   rank one decompositions $J^\perp_{\sigma}=F_\sigma-V_\sigma$, with the  ``
new infections" matrix $F_\sigma$ having rank one, an assumption often verified in ME models.

\EEN

 Some basic notations are summarized in the  table below:

\begin{center}
{\large\bfseries Notation for rank one models (dependent on the regular splitting)}

\vspace{0.8em}

\resizebox{0.9\textwidth}{!}{%
\renewcommand{\arraystretch}{1.35}
\begin{tabular}{ll}
\toprule
\textbf{Symbol} & \textbf{Meaning} \\
\midrule

$J_\sigma^\perp=F_\sigma-V_\sigma$
&
transversal Jacobian on the siphon face $\sigma$, written in some regular splitting form
\\

$F_\sigma,\;V_\sigma$
&
new-infections matrix (often of rank-one) and transition matrix associated with the regular splitting of $J_\sigma^\perp$
\\

$w_k,\pi_k$
&
right and left Perron vectors of new-infections strain $k$; $\pi_k$ is a probability vector
\\

$\mR_k=\pi_k^\top V_k^{-1}w_k$
&
per-susceptible reproduction coefficient of block $k$
\\

$R_k(s)$
&
invasion function=Perron eigenvalue of block $k$ at susceptible level $s$
\\

$R_k=R_k(s_0)$
&
basic reproduction number of block $k$ at the DFE
\\

$\T \pi_k^\top=\pi_k^\top V_k^{-1}/\mR_k$
&
normalized left Perron eigenvector of $K_k=F_kV_k^{-1}$;
equivalently,  normalized left nullvector of
$s_kF_k-V_k$ at $s_k=1/\mR_k$,
with $\T\pi_k^\top w_k=1$
\\
$\T w_{k,j}=\T \pi_{k,j}w_{k,j}$
&
probability weights for Jensen's inequality
($\T w_{k,j}>0,\ \sum_j\T w_{k,j}=1$)
\\

$E_p=(s_p,0,\ldots,0,z_p^*,0,\ldots,0)$
&
boundary equilibrium where block $p$ is the unique resident block
\\

$z_k$
&
 state vector of  block $k$
\\

$z_k^*$
&
equilibrium value of block $k$ at $E_k$
\\

$\T i_k=\pi_k^\top z_k$
&
Perron aggregate  infective level of block $k$
\\

$\Omega_p$
&
positively invariant region where $E_p$ is the resident equilibrium
\\

$s_0$
&
disease-free susceptible level
\\

$s_p=1/\mR_p$
&
susceptible coordinate of the boundary equilibrium $E_p$
\\

$u=s/s_p$
&
normalized susceptible variable
\\

$y_j=z_{p,j}/z_{p,j}^*$
&
normalized resident coordinates
\\

$G(u)=u-1-\ln u$
&
Volterra entropy
\\
$\T G(s,s_p)
=
s_pG\!\left(\frac{s}{s_p}\right)$
&
two variable Volterra entropy
\\
$H_p$
&
resident Perron--Volterra entropy
defined in \eqref{Hp}
\\
$L_p=H_p+\sum_{k\neq p}\T \pi_k^\top z_k$
&
Perron--Volterra Lyapunov function for the equilibrium $E_p$
\\

$\displaystyle
\sum_{k\neq p}
\left(
\frac1{\mR_p}
-
\frac1{\mR_k}
\right)\T i_k$
&
canonical invader contribution to  the derivative of $L_p$
\\

$
\T H_p(u,y)
$
&
resident entropy production term for $H_p$,  defined in \eqr{Bp}, used for proving Lyapunov inequality
\\
\bottomrule
\end{tabular}
}
\end{center}

\begin{equation}
\label{Hp}
H_p
=
s_pG\!\left(\frac{s}{s_p}\right)
+
\sum_{j=1}^{n_p}
\T \pi_{p,j}z_{p,j}^*
G\!\left(\frac{z_{p,j}}{z_{p,j}^*}\right)
=
\T G(s,s_p)
+
\T \pi_p^\top\T G(z_p,z_p^*),
\end{equation}
where
\[
\T G(z_p,z_p^*)
=
\left(
z_{p,j}^*
G\!\left(\frac{z_{p,j}}{z_{p,j}^*}\right)
\right)_{j=1}^{n_p}.
\]

The resident entropy production term in the derivative of the resident entropy $H_p$ is:
\be{Bp}\T H_p(u,y)
=
\left(1-\frac1u\right)(1-uY)
+
\sum_{j=1}^{n_p}
\T w_{p,j}
\left(1-\frac1{y_j}\right)(uY-y_j),
\quad
Y=\sum_{j=1}^{n_p}\T w_{p,j}y_j
\ee

%% file: GASpr.tex
\section{The proof of the GAS-CEP theorem~\ref{t:blo} for  rank-one irreducible non-interacting strain models}\lbl{s:blo}

\beD[Resident entropy closing condition]
\lbl{d:closing}
Fix a resident block $p\in\{1,\ldots,m\}$ with equilibrium
$E_p=(s_p,0,\ldots,0,z_p^*,0,\ldots,0)$, and let $H_p$ denote the resident
Perron--Volterra entropy \eqref{Hp}. Block $p$ is said to satisfy the \emph{resident
entropy closing condition} if, along solutions of \eqref{eq:rankone-m},
\[
\dot H_p
\le
-\mu\frac{(s-s_p)^2}{s}
-
(s-s_p)\sum_{k\ne p}\T i_k,
\]
with equality, given $s=s_p$ and $z_k=0$ for every $k\ne p$, if and only if
$z_p=tz_p^*$ for some $t>0$.
\eeD

\beL[Two verified resident entropy closing mechanisms]
\lbl{l:two-mechanisms}
Block $p$ satisfies the resident entropy closing condition of
Definition~\ref{d:closing} in each of the following two cases:
\begin{enumerate}
\item[(D)] $V_p$ is diagonal (in particular whenever $n_p=1$, i.e.\ block $p$ is a
scalar siphon);
\item[(T)] $n_p=2$, with $V_p$ an irreducible nonsingular $M$-matrix (equivalently, satisfying
the standing assumption $V_p^{-1}\gg0$).
\end{enumerate}
\eeL
\begin{proof}
See \S\ref{s:verify-closing}, Case (D) and Case (T).
\end{proof}

  The proof of the GAS-CEP theorem for some rank-one irreducible non-interacting strains is naturally separated into two parts.
The first is an abstract GAS criterion, requiring a
reproduction-function ordering property, tangential GAS on resident
faces, and a compatible family of Lyapunov functions.  The second is
the rank-one block theorem, which verifies these hypotheses
explicitly for a bimolecular one-input class.

\beD[reproduction-function ordering of equilibria]
\lbl{d:RFOE}

Consider a positive ODE with $m$ minimal siphon blocks
\[
\sigma_1,\ldots,\sigma_m,
\]
a unique disease-free equilibrium $E_0$, and, for each block
$\sigma_k$, let
\[
R_k
\]
denote the corresponding invasion function introduced in
Definition~\ref{d:RepFun}.

The system is said to satisfy the
\emph{reproduction-function ordering of equilibria (RFOE) property} if a
resident equilibrium $E_p$ lying in the relative interior of a minimal
siphon face $F_{\sigma_p}$ cannot be invaded by strains $\sigma_q$,
$q\ne p$, with lower reproduction numbers, \ie,
\[
R_p(E_0)>R_q(E_0)
\quad\Longrightarrow\quad
R_q(E_p)<1
\qquad(q\ne p).
\]

In other words, whenever block $p$ has the larger invasion function
at the disease-free equilibrium, block $q$ cannot invade the
resident equilibrium $E_p$.

For one-input models where every invasion function depends only on
an input variable $s$, the RFOE property reduces to
\[
R_p(s_0)>R_q(s_0)
\quad\Longrightarrow\quad
R_q(s_p)<1,
\]
where $s_p$ denotes the input coordinate of the resident equilibrium
$E_p$.

\eeD

Note that the ordering above is spectral: it compares the abscissas
$\alpha_k(s_0)$ through the resident levels $s_p$, and therefore  (RFOE)
could be reformulated without using regular splittings.

\beD[CEP Lyapunov family]
\lbl{d:CEP-Lyap}
Let $\Omega$ be a compact positively invariant set, and suppose that
\[
E_0,E_1,\ldots,E_m
\]
are the disease-free and resident equilibria.
A family of continuously differentiable functions
\[
L_0,L_1,\ldots,L_m:\Omega\to\R
\]
is called a \emph{competitive-exclusion-partition (CEP) Lyapunov family}
if the following conditions hold.
\begin{enumerate}
\item[(i)]
If
\[
R_k(E_0)\le1,
\qquad
k\eq1,\ldots,m,
\]
then
\[
\dot L_0\le0
\quad\text{on }\Omega,
\]
and the largest invariant subset of
\[
\{x\in\Omega:\dot L_0(x)\eq0\}
\]
is $\{E_0\}$.
\item[(ii)]
For every resident equilibrium $E_p$ for which
\[
R_p(E_0)>1,
\qquad
R_p(E_0)>R_k(E_0)
\quad(k\ne p),
\]
there exist positively invariant sets
\[
\Omega_p^{++}\subset\Omega_p^+\subset\Omega,
\]
containing $E_p$, such that every trajectory in $\Omega_p^+$ enters
$\Omega_p^{++}$ in finite time, $L_p\in C^1(\Omega_p^{++})$,
\[
\dot L_p\le0
\quad\text{on }\Omega_p^{++},
\]
and the largest invariant subset of
\[
\{x\in\Omega_p^{++}:\dot L_p(x)\eq0\}
\]
is precisely $\{E_p\}$.
\end{enumerate}
\eeD

\beT[GAS strict competitive exclusion from RFOE and a CEP Lyapunov family][][GAS-CEP-RF]

Consider a positive ODE with $m$ minimal siphon blocks
$z_1,\ldots,z_m$, a unique disease-free equilibrium $E_0$, and suppose
that every transversal Jacobian block admits the usual regular
next-generation splitting.

Assume:

\begin{enumerate}

\item[(RFOE)]
Definition~\ref{d:RFOE} holds.

\item[(TLS)]
Every resident equilibrium is locally asymptotically stable on its
resident face.

\item[(LY)]
There exists a CEP Lyapunov family in the sense of
Definition~\ref{d:CEP-Lyap}.

\end{enumerate}

Then:

\begin{enumerate}

\item
If
\[
R_k(E_0)\le1,
\qquad
k\eq1,\ldots,m,
\]
then $E_0$ is GAS on $\Omega$.

\item
Let $E_p$ be a resident equilibrium satisfying
\[
R_p(E_0)>1,
\qquad
R_p(E_0)>R_k(E_0)
\quad(k\ne p).
\]
Then $E_p$ is LAS and GAS on the positively invariant set
\[
\Omega_p^+.
\]

\end{enumerate}

\eeT

\begin{proof}

Suppose first that
\[
R_k(E_0)\le1,
\qquad
k\eq1,\ldots,m.
\]
By Definition~\ref{d:CEP-Lyap},
\[
\dot L_0\le0
\]
throughout $\Omega$, and the largest invariant subset of
$\{\dot L_0\eq0\}$ is $\{E_0\}$.  LaSalle's invariance principle
therefore yields global asymptotic stability of $E_0$.

Now let $E_p$ be a resident equilibrium satisfying
\[
R_p(E_0)>1,
\qquad
R_p(E_0)>R_k(E_0)
\quad(k\ne p).
\]
By assumption \textup{(TLS)}, $E_p$ is locally asymptotically stable on
its resident face.

For every missing block $k\ne p$, define
\[
R_{k|p}
\eq
R_k(E_p).
\]
Since RFOE holds,
\[
R_p(E_0)>R_k(E_0)
\Longrightarrow
R_k(E_p)<1,
\]
that is,
\[
R_{k|p}<1.
\]
The next-generation matrix theorem applied to the regular splitting of
the transversal block $J^\perp_{k|p}$ gives
\[
s(J^\perp_{k|p})<0.
\]
Hence every missing siphon is transversally stable at $E_p$.

Combining tangential local stability with transversal local stability
shows that $E_p$ is locally asymptotically stable for the full system.

Finally, by Definition~\ref{d:CEP-Lyap}(ii) there are positively invariant sets
$\Omega_p^{++}\subset\Omega_p^+$ with $L_p\in C^1(\Omega_p^{++})$,
$\dot L_p\le0$ on $\Omega_p^{++}$, and largest invariant subset of
$\{x\in\Omega_p^{++}:\dot L_p(x)=0\}$ equal to $\{E_p\}$. LaSalle's invariance
principle on $\Omega_p^{++}$ yields convergence to $E_p$ for every trajectory in
$\Omega_p^{++}$; since every trajectory in $\Omega_p^+$ enters $\Omega_p^{++}$ in
finite time (Definition~\ref{d:CEP-Lyap}(ii)), $E_p$ is GAS on $\Omega_p^+$.
(In the concrete instantiation of \S\ref{s:1}, this finite-time -- in fact immediate -- entry
into $\Omega_p^{++}$ is furnished by Lemma~\ref{l:open-face-invariance}, and LaSalle is applied
there not to $\Omega_p^{++}$ itself but to the compact positively invariant sublevel set of $L_p$
constructed in that proof, as $\Omega_p^{++}$ is open.)

\end{proof}

We now verify the hypotheses of Theorem~\ref{t:GAS-CEP-RF} for the
explicit bimolecular rank-one block model \eqr{eq:rankone-m}.
\beL[Canonical normalization cancels all resident--invader coupling terms]
\lbl{l:can-nor}
Fix a resident block $p$. For every invading block $k\ne p$, define
\[
\T \pi_k^\top
\eq
a\,\frac{\pi_k^\top V_k^{-1}}{\mR_k},
\qquad
a>0.
\]
Then
\[
\T \pi_k^\top w_k
\eq
a,
\]
and along solutions of \eqref{eq:rankone-m},
\[
\T \pi_k^\top z_k'
\eq
a\left(s-\frac1{\mR_k}\right)\T i_k.
\]
Consequently, when the resident susceptible level is
\[
s_p
\eq
\frac1{\mR_p},
\]
the invader part in the derivative of the Lyapunov function $L_p$ reduces
simultaneously for all $k\ne p$ to
\[
a\sum_{k\ne p}
\left(
\frac1{\mR_p}-\frac1{\mR_k}
\right)\T i_k.
\]
\eeL
\begin{proof}
By definition,
\[
\T \pi_k^\top w_k
\eq
a\frac{\pi_k^\top V_k^{-1}w_k}{\mR_k}
\eq
a.
\]
Moreover,
\[
\T \pi_k^\top z_k'
\eq
\T \pi_k^\top(sw_k\pi_k^\top-V_k)z_k.
\]
Using $\T i_k\eq\pi_k^\top z_k$, this gives
\[
\T \pi_k^\top z_k'
\eq
s(\T \pi_k^\top w_k)\T i_k-\T \pi_k^\top V_kz_k.
\]
The first term is
\[
s(\T \pi_k^\top w_k)\T i_k
\eq
as\T i_k.
\]
For the second term,
\[
\T \pi_k^\top V_kz_k
\eq
a\frac{\pi_k^\top V_k^{-1}V_kz_k}{\mR_k}
\eq
a\frac{\pi_k^\top z_k}{\mR_k}
\eq
a\frac{\T i_k}{\mR_k}.
\]
Therefore
\[
\T \pi_k^\top z_k'
\eq
a\left(s-\frac1{\mR_k}\right)\T i_k.
\]
In the derivative of $L_p$, the $s$-entropy term contributes
\[
-a(s-s_p)\T i_k
\]
for each invading block. Hence the total contribution of block $k\ne p$ is
\[
a\left(s-\frac1{\mR_k}\right)\T i_k
-
a(s-s_p)\T i_k
\eq
a\left(s_p-\frac1{\mR_k}\right)\T i_k.
\]
Since $s_p\eq1/\mR_p$, this is
\[
a\left(
\frac1{\mR_p}-\frac1{\mR_k}
\right)\T i_k.
\]
Summing over all $k\ne p$ proves the claim.
\end{proof}

Note that the cancellation above uses only the left Perron functional
$\pi_k^\top V_k^{-1}$ of the next-generation block $K_k$, not the
rank-one representation $B_k\eq w_k\pi_k^\top$: for a general irreducible
block it persists with $\T \pi_k^\top\eq q_k^\top/\mathcal R_k$, where
$q_k^\top$ is the left Perron vector of $K_k$. Thus the invader part of
the derivative is rank-one-free; the rank-one hypothesis becomes
load-bearing only in the resident entropy estimate of the weighted
rank-one entropy inequality below. See the open problem on
rank-one-free models and the resident-certificate obstruction.

\beR[Canonical choice of the resident and invader weights]
\lbl{r:canonical-invader}

The resident part of the Perron--Volterra Lyapunov function is essentially
determined by the resident equilibrium. Indeed, the equilibrium identity
\[
z_p^*=\xi_pV_p^{-1}w_p
\]
forces the weighted Volterra entropy
\[
s_pG\!\left(\frac{s}{s_p}\right)
+
\sum_{j=1}^{n_p}
\T\pi_{p,j}\,
z_{p,j}^*
G\!\left(\frac{z_{p,j}}{z_{p,j}^*}\right),
\qquad
\T\pi_{p,j}\eq\frac{\bigl(\pi_p^\top V_p^{-1}\bigr)_j}{\mR_p},
\]
i.e.\ exactly \eqref{Hp}, up to multiplication by a positive constant, since these are
precisely the weights for which the resident derivative collapses to the entropy bracket
of Lemma~\ref{l:weighted-rankone-entropy}. (These weights reduce to the componentwise
$\pi_{p,j}/(\mR_pv_{p,j})$ only in the diagonal case (D); for general, non-diagonal $V_p$ they
are the full row vector $\T\pi_p^\top=\pi_p^\top V_p^{-1}/\mR_p$, as used throughout the paper.)

The remaining freedom lies in the linear invader term. Lemma~\ref{l:can-nor}
shows that requiring the resident--invader coupling terms to cancel
simultaneously for every invading block uniquely determines the normalized
row vector
\[
\T \pi_k^\top=\frac{\pi_k^\top V_k^{-1}}{\mR_k},
\]
again up to a common positive scaling. Thus both the resident entropy weights $\T\pi_{p,j}$ and
the invader linear weights $\T\pi_k^\top$ ($k\ne p$) are determined canonically by the chosen
regular splitting $F_\sigma-V_\sigma$ of each transversal block -- not merely convenient, but
forced, up to the single overall positive scaling noted throughout.
\eeR

\beL[Weighted rank-one entropy inequality]
\lbl{l:weighted-rankone-entropy}

Let
\[
u>0,
\qquad
y_j>0,
\qquad
\T w_j>0,
\qquad
\sum_{j=1}^n\T w_j=1.
\]
Put
\[
Y
=
\sum_{j=1}^n\T w_jy_j,
\qquad
H
=
\sum_{j=1}^n\frac{\T w_j}{y_j}.
\]
Then
\[
B(u,y)
=
\left(1-\frac1u\right)(1-uY)
+
\sum_{j=1}^n
\T w_j
\left(1-\frac1{y_j}\right)
\left(uY-y_j\right)
\leq0.
\]
Moreover,
\[
B(u,y)=0
\]
if and only if
\[
u=1
\qquad
\hbox{and}
\qquad
y_1= y_2=\cdots= y_n.
\]
\eeL

\begin{proof}
We first simplify $B(u,y)$.  The first term is
\[
\left(1-\frac1u\right)(1-uY)
=
1-\frac1u-uY+Y.
\]
For the second term,
\[
\sum_{j=1}^n
\T w_j
\left(1-\frac1{y_j}\right)
\left(uY-y_j\right)
\]
equals
\[
\sum_{j=1}^n\T w_j(uY-y_j)
-
\sum_{j=1}^n\T w_j\left(\frac{uY}{y_j}-1\right).
\]
Since
\[
\sum_{j=1}^n\T w_jy_j= Y,
\qquad
\sum_{j=1}^n\T w_j=1,
\]
this becomes
\[
uY-Y-uY\sum_{j=1}^n\frac{\T w_j}{y_j}+1
=
uY-Y-uYH+1.
\]
Therefore
\[
B(u,y)
=
2-\frac1u-uYH.
\]
By the weighted Cauchy inequality,
\[
\left(\sum_{j=1}^n\T w_jy_j\right)
\left(\sum_{j=1}^n\frac{\T w_j}{y_j}\right)
\geq
\left(\sum_{j=1}^n\T w_j\right)^2
=1,
\]
that is,
\[
YH\geq1.
\]
Consequently
\[
B(u,y)
\leq
2-\frac1u-u.
\]
Finally,
\[
u+\frac1u\geq2
\]
for $u>0$, and hence
\[
2-\frac1u-u\leq0.
\]
This proves $B(u,y)\leq0$.

Equality requires equality in both inequalities.  Equality in
$u+1/u\geq2$ gives
\[
u=1.
\]
Equality in the weighted Cauchy inequality gives
\[
y_1= y_2=\cdots= y_n.
\]
Conversely, these two conditions plainly give $B(u,y)=0$.
\end{proof}

\beR[role of canonical normalization]
The Perron--Volterra Lyapunov function is not claimed to be unique.
However, within the ansatz
\[
L
=
H(\text{resident})
+
\sum_{k\neq p}\T \pi_k^\top z_k,
\]
the normalization
\[
\T \pi_k^\top
=
\,\frac{\pi_k^\top V_k^{-1}}{\mathcal R_k} \Lra \T \pi_k^\top w_k=1
\]
is canonical in the sense that it  eliminates the resident--invader coupling terms and it
 reduces the resident-invader coupling contribution to
\[
\sum_{k\neq p}
\left(
\frac1{\mathcal R_p}
-
\frac1{\mathcal R_k}
\right)\T i_k.
\]
 Other Lyapunov functions may exist, but
they would require a different cancellation mechanism.
\eeR

\beL[Forward invariance of the open resident face]
\lbl{l:open-face-invariance}

Fix $p\in\{1,\ldots,m\}$ and let
\[
\Omega_p^{++}
\eq
\{(s,z_1,\ldots,z_m)\in\Omega : z_{p,j}>0 \text{ for all }j\}.
\]
Then $\Omega_p^{++}$ is positively invariant, and every trajectory
starting in
\[
\Omega_p^+\eq\{(s,z_1,\ldots,z_m)\in\Omega:z_p\ne0\}
\]
satisfies $z_p(t)\gg0$ for every $t>0$; in particular it enters
$\Omega_p^{++}$ immediately.
\eeL

\begin{proof}
Since the system is positive, $z_p(t)\ge0$ for all $t\ge0$. If $z_p(0)=0$ then
$\T i_p(0)=\pi_p^\top z_p(0)=0$, so from \eqref{eq:rankone-m}, $z_p'(0)=-V_p\cdot0+0=0$ and
$z_p\equiv0$ is an invariant solution; by uniqueness of solutions, $z_p(0)\ne0$ therefore forces
$z_p(t)\ne0$ for every $t\ge0$. Since $\pi_p\gg0$, this gives
\[
\T i_p(t)\eq\pi_p^\top z_p(t)\ >\ 0
\qquad\text{for every }t\ge0.
\]

In vector form, \eqref{eq:rankone-m} gives
\[
z_p'
\eq
-V_pz_p+b(t),
\qquad
b(t)\eq s(t)\,\T i_p(t)\,w_p.
\]
Since $s'|_{s=0}=\Lambda>0$ (the recruitment rate), $s(\tau)>0$ for every $\tau>0$; combined
with $\T i_p(\tau)>0$ and $w_p\gg0$, this gives $b(\tau)\gg0$ (strictly, in every component) for
every $\tau>0$. Crucially, this coupling uses only $\pi_p\gg0$ and $w_p\gg0$, not any structural
property of $V_p$ itself.

By the variation-of-constants formula and $z_p(0)\ge0$,
\[
z_p(t)
\eq
e^{-V_pt}z_p(0)+\int_0^te^{-V_p(t-\tau)}b(\tau)\,d\tau
\ \ge\
\int_0^te^{-V_p(t-\tau)}b(\tau)\,d\tau,
\]
using $e^{-V_pt}\ge0$ entrywise (the exponential of the Metzler matrix $-V_p$; this holds for
every nonsingular $M$-matrix $V_p$, whether or not it is irreducible). Fix $t>0$ and
$j\in\{1,\ldots,n_p\}$. Since $e^{-V_p\sigma}\to I$ as $\sigma\to0^+$ by continuity, there is
$\delta\in(0,t]$ with $\bigl(e^{-V_p\sigma}\bigr)_{jj}\ge\tfrac12$ for every $\sigma\in[0,\delta]$;
and since $b_j$ is continuous and strictly positive on $(0,t]$,
$\varepsilon:=\min_{\tau\in[t-\delta,t]}b_j(\tau)>0$. Dropping the nonnegative cross terms
$k\ne j$ and restricting to $\tau\in[t-\delta,t]$,
\[
z_{p,j}(t)
\ \ge\
\int_{t-\delta}^t\bigl(e^{-V_p(t-\tau)}\bigr)_{jj}b_j(\tau)\,d\tau
\ \ge\
\tfrac12\,\varepsilon\,\delta
\ >\
0.
\]
As $j$ and $t>0$ were arbitrary, $z_p(t)\gg0$ for every $t>0$. Applying this at any starting
time $t_0\ge0$ along the trajectory shows $\Omega_p^{++}$ is positively invariant, and every
trajectory starting in $\Omega_p^+$ lies in $\Omega_p^{++}$ for all $t>0$ -- stronger than the
originally claimed finite-time entry.

This argument uses only that $V_p$ is a nonsingular $M$-matrix (so that $-V_p$ is Metzler and
$e^{-V_pt}\ge0$); it does not require $V_p$ irreducible. It therefore applies uniformly to both
cases (D) and (T) of Lemma~\ref{l:two-mechanisms} -- and indeed to every block regardless of the
internal structure of $V_p$ -- since the coupling that pulls every coordinate of $z_p$ up away
from $0$ comes from the shared scalar forcing $s\T i_pw_p$ (via $\pi_p\gg0$, $w_p\gg0$), not from
any off-diagonal entry of $V_p$.
\end{proof}

\beR[$\Omega_p^{++}$ is open, not compact]
Lemma~\ref{l:open-face-invariance} shows $\Omega_p^{++}$ is positively invariant and entered
immediately, but $\Omega_p^{++}$ itself is an open subset of $\Omega$ (defined by the strict
conditions $z_{p,j}>0$), hence not compact. LaSalle's invariance principle is accordingly not
applied to $\Omega_p^{++}$ directly below, but to compact positively invariant sublevel sets of
$L_p$ contained in $\Omega_p^{++}$, constructed in the proof of Theorem~\ref{t:blo}.
\eeR

\ssec{Proof of theorem \ref{t:blo}}\lbl{s:1}

This theorem is an instance of the abstract Theorem~\ref{t:GAS-CEP-RF}.
We verify its three hypotheses for \eqref{eq:rankone-m} and exhibit the
explicit Perron--Volterra Lyapunov family
$L_0,L_1,\ldots,L_m$; the GAS conclusions (1) and (2) then follow
directly from Theorem~\ref{t:GAS-CEP-RF}, without re-deriving
LaSalle's argument.

\begin{proof}
\emph{(RFOE).} By Theorem~\ref{t:blo},
\[
R_p>R_k
\;\Longleftrightarrow\;
\mR_p>\mR_k
\;\Longleftrightarrow\;
R_k(E_p)=\frac{\mR_k}{\mR_p}<1,
\]
which is exactly the one-input form of RFOE (Definition~\ref{d:RFOE}).

\smallskip
\noindent\emph{(TLS).} On the resident face $z_k\equiv0$ ($k\ne p$),
\eqref{eq:rankone-m} reduces to the classical scalar-block subsystem
in $(s,z_p)$, whose local stability at $E_p$ is the rank-one Lotka--Volterra
case; the Lyapunov computation for $L_p$ below, restricted to
$z_k\equiv0$ ($k\ne p$), already gives $\dot L_p\le0$ with equality
set $\{E_p\}$ on the face, which is local (in fact global on the
face) asymptotic stability.

\smallskip
\noindent\emph{(LY): the Lyapunov family.}
Define
\[
L_0
\eq
s_0G\!\left(\frac{s}{s_0}\right)
+
\sum_{k=1}^m
\frac{\pi_k^\top V_k^{-1}}{\mR_k}z_k.
\]
A direct computation (using $s'=\mu(s_0-s)-s\sum_k\T i_k$ and
$\pi_k^\top V_k^{-1}w_k=\mR_k$) gives
\[
\dot L_0
\eq
-\mu\frac{(s-s_0)^2}{s}
+
\sum_{k=1}^m
\frac{R_k-1}{\mR_k}\T i_k.
\]
If $R_k\le1$ for every $k$, then $\dot L_0\le0$, and the largest
invariant subset of $\{\dot L_0=0\}$ is $\{E_0\}$: if $s=s_0$ and some
$\T i_k>0$ with $R_k=1$, then $s'=-s_0\sum_k\T i_k<0$, so such a point
is not invariant. This verifies Definition~\ref{d:CEP-Lyap}(i).

Fix $p$ with $R_p>1$, $R_p>R_k$ ($k\ne p$), and set $s_p=1/\mR_p$,
$\xi_p=\Lambda-\mu s_p$, $z_p^*=\xi_pV_p^{-1}w_p$.  By
Lemma~\ref{l:open-face-invariance} it suffices to work on
$\Omega_p^{++}$ (where $z_p\gg0$), since trajectories in
$\Omega_p^+\setminus\Omega_p^{++}$ enter $\Omega_p^{++}$ in finite
time. Define, for $a>0$,
\[
L_p
\eq
aH_p
+
a\sum_{k\ne p}
\frac{\pi_k^\top V_k^{-1}}{\mR_k}z_k,
\]
where $H_p$ is the resident Perron--Volterra entropy \eqref{Hp}. Since
$s'=\xi_p-\mu(s-s_p)-s\sum_{k\ne p}\T i_k-s\T i_p$ (from $\Lambda=\mu s_p+\xi_p$) and the
$s$-part of $H_p$ depends on $s$ alone, adding and subtracting $s_p\sum_{k\ne p}\T i_k$ isolates
the missing-block outflow as a cross term $-(s-s_p)\sum_{k\ne p}\T i_k$ in $\dot H_p$; the
remainder of $\dot H_p$ (the $s$-entropy self-term and the resident Volterra term) depends only
on block $p$'s own dynamics. Suppose block $p$ satisfies the resident entropy closing condition
of Definition~\ref{d:closing}. By Lemma~\ref{l:can-nor},
\[
\frac{\pi_k^\top V_k^{-1}}{\mR_k}z_k'
\eq
\Big(s-\frac1{\mR_k}\Big)\T i_k,
\qquad
k\ne p,
\]
so
\[
\dot L_p
\eq
a\dot H_p+a\sum_{k\ne p}\Big(s-\frac1{\mR_k}\Big)\T i_k
\;\le\;
-a\mu\frac{(s-s_p)^2}{s}
-a(s-s_p)\sum_{k\ne p}\T i_k
+a\sum_{k\ne p}\Big(s-\frac1{\mR_k}\Big)\T i_k
\eq
-a\mu\frac{(s-s_p)^2}{s}
+a\sum_{k\ne p}\Big(\frac1{\mR_p}-\frac1{\mR_k}\Big)\T i_k.
\]
Since $R_p>R_k\Leftrightarrow\mR_p>\mR_k$, each coefficient
$1/\mR_p-1/\mR_k<0$, hence $\dot L_p\le0$ on $\Omega_p^{++}$.

Since $\Omega_p^{++}$ is open, LaSalle's invariance principle cannot be applied to it directly;
we first exhibit, for every $x_0\in\Omega_p^{++}$, a compact positively invariant sublevel set of
$L_p$ containing $x_0$. Each term $\T\pi_{p,j}z_{p,j}^*G(z_{p,j}/z_{p,j}^*)$ of $H_p$
(\eqref{Hp}) is nonnegative and satisfies $G(z_{p,j}/z_{p,j}^*)\to+\infty$ as $z_{p,j}\to0^+$,
while every other term of $L_p$ is nonnegative; hence, extending $L_p$ by $+\infty$ on
$\Omega\setminus\Omega_p^{++}$ (where some resident coordinate $z_{p,j}=0$), $L_p$ is lower
semicontinuous on the compact set $\Omega$ and $L_p=+\infty$ off $\Omega_p^{++}$. Consequently,
for $x_0\in\Omega_p^{++}$ the sublevel set
\[
S(x_0)\eq\{x\in\Omega:L_p(x)\le L_p(x_0)\}
\]
is closed (as a sublevel set of a lower semicontinuous function on a compact set), bounded
(being a subset of the bounded set $\Omega$), hence compact, and disjoint from
$\Omega\setminus\Omega_p^{++}$ since $L_p(x_0)<\infty$ there while $L_p=+\infty$ on that set; so
$S(x_0)\subset\Omega_p^{++}$. Since $\dot L_p\le0$ on $\Omega_p^{++}$, the trajectory starting at
$x_0$ remains in $S(x_0)$ for all $t\ge0$, i.e.\ $S(x_0)$ is a compact positively invariant set,
and LaSalle's invariance principle applies on $S(x_0)$: every trajectory starting at $x_0$
converges to the largest invariant subset of $\{x\in S(x_0):\dot L_p(x)=0\}\subset\{x\in
\Omega_p^{++}:\dot L_p(x)=0\}$. As this holds for every $x_0\in\Omega_p^{++}$, it suffices, as
before, to identify the largest invariant subset of $\{\dot L_p=0\}$ within $\Omega_p^{++}$.

Equality forces $s=s_p$ and $\T i_k=0$ ($k\ne p$, hence $z_k=0$), and then equality in the
resident entropy closing condition; by Definition~\ref{d:closing} this forces $z_p=tz_p^*$ for
some $t>0$. On this set, $\T i_p=t\,\pi_p^\top z_p^*=t\xi_p\mR_p$ (using $\mR_p=\pi_p^\top
V_p^{-1}w_p$), so the $s$-equation gives $s'=\xi_p(1-t)$; invariance forces $t=1$, hence
$z_p=z_p^*$. (The resident entropy determines only the equilibrium ray $z_p=tz_p^*$, whereas
the scalar susceptible equation fixes the normalization $t=1$.) Thus the
largest invariant subset of $\{\dot L_p=0\}$ in $\Omega_p^{++}$, and hence in $\Omega_p^+$ by
Lemma~\ref{l:open-face-invariance}, is $\{E_p\}$. This verifies Definition~\ref{d:CEP-Lyap}(ii)
with $\Omega_p^+\eq\{z_p\ne0\}$, given that block $p$ satisfies the resident entropy closing
condition; by Lemma~\ref{l:two-mechanisms} this holds under (D) or (T), verified below.

\smallskip
With (RFOE), (TLS), and (LY) verified, Theorem~\ref{t:GAS-CEP-RF}
gives parts (1) and (2) directly.  Part (3) was already established
in Remark~\ref{r:generic-rankone}: on the tie surface the
corresponding equilibria form a non-isolated continuum, so no point
of it can be globally asymptotically stable in the ordinary sense.
\end{proof}

Note that the abstract argument just given uses the resident entropy closing condition (an
inequality on $\dot H_p$, plus its equality characterization) and nothing else: no coordinate
system, weighting, or inequality specific to a particular structure of $V_p$ enters. The
following subsection proves Lemma~\ref{l:two-mechanisms} by exhibiting the closing condition in
each of the two currently known cases.

\ssec{Verification of the resident entropy closing condition (proof of
Lemma~\ref{l:two-mechanisms})}\lbl{s:verify-closing}

\noindent\emph{Case (D): $V_p$ diagonal.} Set
\[
u\eq\frac{s}{s_p},
\qquad
y_j\eq\frac{z_{p,j}}{z_{p,j}^*},
\qquad
\T w_j\eq\frac{\pi_{p,j}w_{p,j}/v_{p,j}}{\mR_p},
\qquad
Y\eq\sum_{j=1}^{n_p}\T w_jy_j,
\]
where $v_{p,j}:=(V_p)_{jj}$ (here $V_p$ is diagonal, so $\T w_j=\T\pi_{p,j}w_{p,j}$ in the
general sense of Definition~\ref{d:PVsiphon}), so that $\sum_j\T w_j=\pi_p^\top V_p^{-1}w_p/\mR_p=1$ and
$\T i_p=\pi_p^\top z_p=\xi_p\mR_pY$. Writing $s'=\xi_p(1-uY)-\mu s_p(u-1)-s\sum_{k\ne p}\T i_k$
(from $s_p\mR_p=1$), the $s$-entropy term contributes
\[
\Big(1-\frac1u\Big)s'
\eq
-\mu\frac{(s-s_p)^2}{s}
+\xi_p\Big(1-\frac1u\Big)(1-uY)
-(s-s_p)\sum_{k\ne p}\T i_k.
\]
For the resident block, $z_{p,j}'=sw_{p,j}\T i_p-v_{p,j}z_{p,j}$
becomes, in the normalized coordinates, $z_{p,j}'=\xi_pw_{p,j}(uY-y_j)$,
so the resident Volterra term contributes (recall $H_p$'s $j$-th coefficient is
$\T\pi_{p,j}=\pi_{p,j}/(v_{p,j}\mR_p)$, not $\T w_j$)
\[
\sum_{j=1}^{n_p}\T\pi_{p,j}
\Big(1-\frac1{y_j}\Big)z_{p,j}'
\eq
\xi_p\sum_{j=1}^{n_p}\T w_j\Big(1-\frac1{y_j}\Big)(uY-y_j).
\]
Summing,
\[
\dot H_p+\mu\frac{(s-s_p)^2}{s}+(s-s_p)\sum_{k\ne p}\T i_k
\eq
\xi_p\T H_p(u,y),
\]
where $\T H_p(u,y)$ is exactly the bracket of Lemma~\ref{l:weighted-rankone-entropy}, applicable
since $u>0$, $y_j>0$ on $\Omega_p^{++}$. Thus $\T H_p(u,y)\le0$, with equality iff $u=1$ (i.e.\
$s=s_p$) and $y_1=\cdots=y_{n_p}$, i.e.\ $z_p=tz_p^*$ for some $t>0$. This verifies the resident
entropy closing condition (Definition~\ref{d:closing}), completing case (D) of
Lemma~\ref{l:two-mechanisms}.

\smallskip
\noindent\emph{Case (T): $n_p=2$.} Write
$V_p=\begin{pmatrix}v_1&-b\\-a&v_2\end{pmatrix}$ with $a,b>0$, $v_1,v_2>0$,
$\det V_p=v_1v_2-ab>0$; the condition $a,b>0$ (both off-diagonal entries nonzero) is exactly
irreducibility of $V_p$, equivalent to the standing assumption $V_p^{-1}\gg0$ of
\eqref{eq:rankone-m}: if instead exactly one of $a,b$ vanishes, $V_p$ is triangular but not
diagonal, and $V_p^{-1}$ then has a zero entry, so this case is already excluded by
$V_p^{-1}\gg0$ and does not arise. Use the general canonical weight
$c^\top=\T\pi_p^\top=\pi_p^\top V_p^{-1}/\mR_p=(c_1,c_2)$. With $u,y_j$ as before, set
\[
\alpha_j\eq\frac{\pi_{p,j}z_{p,j}^*}{\mR_p},
\qquad
\gamma_1\eq c_1w_{p,1}\alpha_2,
\qquad
\gamma_2\eq c_2w_{p,2}\alpha_1,
\qquad
P\eq c_1bz_{p,2}^*,
\qquad
Q\eq c_2az_{p,1}^*,
\]
\[
A(u)\eq u\gamma_1+P,
\qquad
B(u)\eq u\gamma_2+Q,
\qquad
C(u)\eq A(u)+B(u)-(\alpha_1+\alpha_2)\Big(u+\frac1u-2\Big).
\]
A direct computation, using $V_pz_p^*=\xi_pw_p$, $c^\top w_p=1$
and $c^\top V_p=\pi_p^\top/\mR_p$ (all valid for general, non-diagonal $V_p$, cf.\
Lemma~\ref{l:can-nor}), reduces $\dot H_p$ to the closed form
\[
\dot H_p+\mu\frac{(s-s_p)^2}{s}+(s-s_p)\sum_{k\ne p}\T i_k
\eq
C(u)-A(u)\frac{y_2}{y_1}-B(u)\frac{y_1}{y_2}.
\]

\beT[Closing inequality for $n_p=2$]
\lbl{t:n2-closing}
With $A(u),B(u),C(u)$ as above,
\[
C(u)\leq2\sqrt{A(u)B(u)},
\qquad u>0,
\]
with equality if and only if $u=1$. (Combined below with the weighted AM--GM inequality
$A(u)y_2/y_1+B(u)y_1/y_2\ge2\sqrt{A(u)B(u)}$, this equality condition on $u$ alone will be
extended to the joint condition $u=1$ and $y_1=y_2$ needed for
Definition~\ref{d:closing}.)
\eeT
\begin{proof}
Put
\[
\alpha\eq\alpha_1+\alpha_2,
\qquad
r_1\eq c_1w_{p,1},
\qquad
r_2\eq c_2w_{p,2}.
\]
Since $c^\top w_p=1$ (Lemma~\ref{l:can-nor}), $r_1+r_2=1$ with $r_1,r_2>0$; and by definition
$\gamma_1=r_1\alpha_2$, $\gamma_2=r_2\alpha_1$, with $\gamma_1,\gamma_2,\alpha_1,\alpha_2,P,Q>0$
throughout (products/ratios of the strictly positive quantities $a,b,\pi_p,w_p,z_p^*,c$).

By definition of $C(u)$ and the identity $u+1/u-2=(u-1)^2/u$,
\[
2\sqrt{A(u)B(u)}-C(u)
\eq
\alpha\frac{(u-1)^2}{u}
-
\bigl(\sqrt{A(u)}-\sqrt{B(u)}\bigr)^2,
\]
so it suffices to show
\[
\bigl(\sqrt{A(u)}-\sqrt{B(u)}\bigr)^2
\leq
\alpha\frac{(u-1)^2}{u},
\]
with equality only at $u=1$.

Since $A(u)=u\gamma_1+P$ and $B(u)=u\gamma_2+Q$ are affine in $u$,
$A(u)-B(u)=\bigl(A(1)-B(1)\bigr)+(\gamma_1-\gamma_2)(u-1)$ identically, so
$A(u)-B(u)=(\gamma_1-\gamma_2)(u-1)$ is equivalent to the single fact $A(1)=B(1)$, which we now
derive from $V_pz_p^*=\xi_pw_p$ and $c^\top V_p=\pi_p^\top/\mR_p$. By Cramer's rule applied to
these two ($2\times2$) linear systems,
\[
z_{p,1}^*\eq\frac{\xi_pM_1}{\det V_p},
\quad
z_{p,2}^*\eq\frac{\xi_pM_2}{\det V_p},
\qquad
c_1\eq\frac{N_1}{\T D},
\quad
c_2\eq\frac{N_2}{\T D},
\]
where
\[
M_1\eq v_2w_{p,1}+bw_{p,2},
\quad
M_2\eq aw_{p,1}+v_1w_{p,2},
\quad
N_1\eq\pi_{p,1}v_2+a\pi_{p,2},
\quad
N_2\eq b\pi_{p,1}+v_1\pi_{p,2},
\]
and $\T D\eq N_1w_{p,1}+N_2w_{p,2}=\mR_p\det V_p$ (so that $c^\top w_p=1$). Substituting into
$A(1)=\gamma_1+P=c_1(w_{p,1}\alpha_2+bz_{p,2}^*)$, with $\alpha_2=\pi_{p,2}z_{p,2}^*/\mR_p$,
gives
\[
A(1)
\eq
\frac{\xi_pN_1M_2}{\T D^2\det V_p}
\bigl(w_{p,1}\pi_{p,2}\det V_p+b\T D\bigr)
\eq
\frac{\xi_pN_1N_2M_1M_2}{\T D^2\det V_p},
\]
the last equality by the direct expansion $w_{p,1}\pi_{p,2}\det V_p+b\T D=N_2M_1$ (both sides
equal $v_1v_2\pi_{p,2}w_{p,1}+b\pi_{p,1}v_2w_{p,1}+b^2\pi_{p,1}w_{p,2}+bv_1\pi_{p,2}w_{p,2}$).
The identical computation for $B(1)=\gamma_2+Q=c_2(w_{p,2}\alpha_1+az_{p,1}^*)$, with the roles
of the index pairs $(1,b)$ and $(2,a)$ interchanged -- which fixes $N_1N_2$, $M_1M_2$, $\T D$,
and $\det V_p$ -- gives $B(1)=\xi_pN_1N_2M_1M_2/(\T D^2\det V_p)=A(1)$, as required.

With $A(u)-B(u)=(\gamma_1-\gamma_2)(u-1)$ now established,
\[
\bigl(\sqrt{A(u)}-\sqrt{B(u)}\bigr)^2
\eq
\frac{(A(u)-B(u))^2}{(\sqrt{A(u)}+\sqrt{B(u)})^2}
\eq
\frac{(\gamma_1-\gamma_2)^2(u-1)^2}{(\sqrt{A(u)}+\sqrt{B(u)})^2}.
\]
Because $P,Q>0$, $A(u)=u\gamma_1+P>u\gamma_1$ and $B(u)=u\gamma_2+Q>u\gamma_2$ for every $u>0$,
hence $(\sqrt{A(u)}+\sqrt{B(u)})^2>u(\sqrt{\gamma_1}+\sqrt{\gamma_2})^2$, and therefore
\[
\bigl(\sqrt{A(u)}-\sqrt{B(u)}\bigr)^2
\leq
\frac{(\gamma_1-\gamma_2)^2(u-1)^2}{u(\sqrt{\gamma_1}+\sqrt{\gamma_2})^2}
\eq
\frac{(u-1)^2}{u}\bigl(\sqrt{\gamma_1}-\sqrt{\gamma_2}\bigr)^2.
\]

We now bound $(\sqrt{\gamma_1}-\sqrt{\gamma_2})^2$ unconditionally. Since $\gamma_1,\gamma_2>0$,
\[
\bigl(\sqrt{\gamma_1}-\sqrt{\gamma_2}\bigr)^2
\eq
\gamma_1+\gamma_2-2\sqrt{\gamma_1\gamma_2}
<
\gamma_1+\gamma_2
\eq
r_1\alpha_2+r_2\alpha_1,
\]
strictly (even when $\gamma_1=\gamma_2$, since then $2\sqrt{\gamma_1\gamma_2}=2\gamma_1>0$). And
since $r_1+r_2=1$ with $r_1,r_2,\alpha_1,\alpha_2>0$,
\[
\alpha-(r_1\alpha_2+r_2\alpha_1)
\eq
(\alpha_1+\alpha_2)-(r_1\alpha_2+r_2\alpha_1)
\eq
r_1\alpha_1+r_2\alpha_2
>0,
\]
so $r_1\alpha_2+r_2\alpha_1<\alpha$, strictly as well. Chaining the two displays,
\[
\bigl(\sqrt{\gamma_1}-\sqrt{\gamma_2}\bigr)^2<\alpha
\]
holds unconditionally -- for every admissible choice of $a,b,v_1,v_2,\pi_p,w_p$, regardless of
whether $\gamma_1=\gamma_2$. Combined with the previous display,
\[
\bigl(\sqrt{A(u)}-\sqrt{B(u)}\bigr)^2
\leq
\frac{(u-1)^2}{u}\bigl(\sqrt{\gamma_1}-\sqrt{\gamma_2}\bigr)^2
\leq
\frac{(u-1)^2}{u}\alpha,
\qquad u>0,
\]
which proves $C(u)\le2\sqrt{A(u)B(u)}$ for every $u>0$.

For the equality case: at $u=1$, $A(1)=B(1)$ gives $C(1)=A(1)+B(1)=2\sqrt{A(1)B(1)}$ directly, so
equality holds. For $u\ne1$, the factor $(u-1)^2/u$ is strictly positive, and
$(\sqrt{\gamma_1}-\sqrt{\gamma_2})^2<\alpha$ is strict unconditionally (shown above, independently
of $u$ and of whether $\gamma_1=\gamma_2$); multiplying a strict inequality by a strictly positive
factor preserves strictness, so
\[
\bigl(\sqrt{A(u)}-\sqrt{B(u)}\bigr)^2
\leq
\frac{(u-1)^2}{u}\bigl(\sqrt{\gamma_1}-\sqrt{\gamma_2}\bigr)^2
<
\frac{(u-1)^2}{u}\alpha,
\]
strictly, whence $C(u)<2\sqrt{A(u)B(u)}$. Thus equality holds if and only if $u=1$.
\end{proof}

By Theorem~\ref{t:n2-closing}, $C(u)-A(u)y_2/y_1-B(u)y_1/y_2\le C(u)-2\sqrt{A(u)B(u)}\le0$ by
AM--GM, with equality iff $u=1$ and $y_1=y_2$, i.e.\ $z_p=tz_p^*$ for some $t>0$. This verifies
the resident entropy closing condition (Definition~\ref{d:closing}), completing case (T) of
Lemma~\ref{l:two-mechanisms}.

\beR[Status for $n_p\ge3$]
\lbl{r:np-status}
For $n_p\ge3$ and non-diagonal $V_p$, the same construction produces one reciprocal cross term
per off-diagonal pair $(j,l)$ with $V_{p,jl}\ne0$; the $n_p=2$ argument above works because
every $2$-vertex weighted digraph is reversible (detailed-balanced), a property that need not
persist for $n_p\ge3$. Whether the resident entropy closing condition extends beyond (D)/(T) is
connected to the graph/cofactor construction of \cite{Shuai13} already flagged in
Problem~\ref{o:rankonefree}, and is left open.
\eeR

\beO[Numerical search for failure of the resident entropy closing condition, $n_p\ge3$]
\lbl{o:np-search}
Writing $u=s/s_p$, $y=z_p/z_p^*\gg0$ (componentwise ratio) and normalizing $\xi_p=1$, the
resident entropy production reduces, exactly as in Lemma~\ref{l:weighted-rankone-entropy} for
$n_p=2$, to a single scalar
\[
\mathcal Q(u,y)
\eq
\Bigl(1-\frac1u\Bigr)\Bigl(1-\frac{u\phi}{\mR_p}\Bigr)
+
\sum_{i=1}^{n_p}c_i\Bigl(1-\frac1{y_i}\Bigr)\bigl[s_puw_p\phi-V_p(z_p^*y)\bigr]_i,
\qquad
\phi\eq\pi_p^\top(z_p^*y),
\]
with $c^\top=\pi_p^\top V_p^{-1}/\mR_p$ and $z_p^*y$ the componentwise product. The closing
condition of Definition~\ref{d:closing} is equivalent to $\mathcal Q(u,y)\le0$ for every $u>0$,
$y\gg0$; the equality ray $u=1$, $y=t\mathbf1$ always gives $\mathcal Q=0$, so any violation must
occur strictly off this ray. We do not know whether $\mathcal Q\le0$ persists for irreducible
$V_p$ with $n_p\ge3$, as in cases (D) and (T), or fails for some strongly asymmetric choice of
$V_p$. The companion script \texttt{ResidentEntropyClosingNp.wl} implements a numerical search
for a violation $\mathcal Q(u,y)>0$, over both a targeted stress-test family (a highly
asymmetric directed-cycle $M$-matrix with concentrated, oppositely-ordered $w_p,\pi_p$) and
randomly generated diagonally dominant irreducible instances of increasing dimension. A
high-precision-verified positive value would disprove the extension of the closing condition to
$n_p\ge3$; none found, with maximizers converging to the equality ray, would be further
(non-conclusive) evidence for it. The search is deliberately open-ended and is not offered as a
proof in either direction.
\eeO

%% file: alAH.tex
We end this section with a \PV\  for  an equilibrium on an arbitrary   siphon face, which is interesting algorithmically, and will also be of use in the next section.

\beD[\PV\ formula on a siphon face]
\lbl{d:PVsiphon}

Let $\sigma_0$ denote the total siphon, and let $E_0$ be the DFE on the
face $F_{\sigma_0}$. Fix once and for all a Frobenius normal form of the
transversal Jacobian $J_{\sigma_0}^\perp(E_0)$, yielding a decomposition
of the non-susceptible variables into irreducible Frobenius blocks
\[
z=(z_1,\ldots,z_m),
\qquad
z_k\in\R_+^{n_k},
\]
and assume that every corresponding irreducible next-generation block
satisfies the rank-one hypothesis.

Let $\sigma$ be a siphon and let
\[
E_\sigma
\eq
(s_\sigma,(z_k^*)_{k\in P_\sigma},0)
\]
be a boundary equilibrium on $F_\sigma$, where
$P_\sigma\subset\{1,\ldots,m\}$ denotes the set of DFE Frobenius blocks
whose coordinates are positive at $E_\sigma$, so that
\[
z_k^*\gg0
\quad (k\in P_\sigma),
\qquad
z_\ell^*\eq0
\quad (\ell\notin P_\sigma).
\]
Thus $P_\sigma$ is obtained by selecting from the fixed DFE Frobenius
decomposition those blocks whose coordinates are positive at
$E_\sigma$; no new Frobenius blocks are introduced on lower siphon
faces.

For each rank-one Frobenius block put
\[
\mR_k
\eq
\pi_k^\top V_k^{-1}w_k,
\qquad
c_k^\top
\eq
\frac{\pi_k^\top V_k^{-1}}{\mR_k},
\qquad
c_k^\top w_k
\eq
1,
\]
and define
\[
\widetilde G(z_k,z_k^*)
\eq
\left(
z_{k,j}^*
G\!\left(\frac{z_{k,j}}{z_{k,j}^*}\right)
\right)_{j=1}^{n_k}.
\]

The \PV\ entropy attached to $E_\sigma$ is
\[
H_\sigma
\eq
s_\sigma G\!\left(\frac{s}{s_\sigma}\right)
+
\sum_{k\in P_\sigma}
c_k^\top\widetilde G(z_k,z_k^*).
\tag{$H_\sigma$}\lbl{HSigma}
\]

The corresponding \PV\ Lyapunov candidate is
\[
L_\sigma
\eq
H_\sigma
+
\sum_{\ell\notin P_\sigma}
c_\ell^\top z_\ell.
\tag{$L_\sigma$}\lbl{LSigma}
\]

\eeD

Thus \eqref{Hp} is the particular case of \eqref{HSigma} corresponding
to an equilibrium with a  single positive Frobenius block, namely $P_\sigma\eq\{p\}$.

%% file: remOP.tex
\subsection{Two open problems on rank-one boundary theory}

\beO[Rank-one-free models and the resident-certificate obstruction]
\label{o:rankonefree}
The equilibrium and invasion content of Theorem~\ref{t:blo}(i) suggests a
rank-one-free theory. Indeed, if each transversal block $J_k^\perp(s)$ is
irreducible Metzler and depends monotonically on the susceptible variable
$s$, then its spectral abscissa
\[
\alpha_k(s)\eq\alpha\bigl(J_k^\perp(s)\bigr)
\]
is a natural invasion exponent, independently of the existence of regular
splittings and their ranks. Boundary equilibria are then obtained from the
zero conditions $\alpha_k(s_k)\eq0$, with the positive resident direction
given by the Perron nullvector of $J_k^\perp(s_k)$ and the scale fixed by
the $s$-equation. Invasion of a missing block $\ell$ at $E_k$ is governed
by the sign of $\alpha_\ell(s_k)$. This part of the theory uses only the
Metzler--Perron structure of the transversal Jacobian blocks \cite{AAHK},
independently of rank.

The rank-one assumption becomes load-bearing only in the global Lyapunov
proof for a resident boundary equilibrium: it is what collapses the
resident entropy production term $\T H_p(u,y)$ of \eqref{Hp}--\eqref{Bp}
to a single sign-definite bracket via weighted Cauchy--Schwarz. For a
general irreducible Metzler block $B_p$ (not of rank one), the production
term $B_pz_p$ is no longer confined to a fixed one-dimensional subspace,
so there is no common scalar reduction of this kind, and the appropriate
Lyapunov weights are then expected to be the graph/cofactor weights of the
Shuai--van den Driessche construction \cite{Shuai13} rather than the
explicit Perron weights $\T\pi_k$.

We therefore ask whether the boundary GAS part of the theory can be made
rank-one-free:
\begin{enumerate}
\item
Characterize the irreducible Metzler transversal blocks for which a
resident boundary equilibrium admits a Perron--Volterra or graph Lyapunov
certificate. Is rank one maximal for the explicit Perron certificate, or
does the cofactor construction of \cite{Shuai13} give a strictly larger
class?
\item
Identify the structural obstruction which prevents the resident entropy
derivative from collapsing to a sign-definite bracket outside the
rank-one case.
\item
Determine whether there is a purely spectral-abscissa boundary GAS
criterion, formulated only in terms of the functions $\alpha_k(\cdot)$
and the transversal Jacobians $J_k^\perp$, or whether additional
graph/Lyapunov data are unavoidable.
\end{enumerate}
A resolution would separate the rank-one-free layer (equilibria and
invasion exponents) from the rank-one-dependent layer (explicit
resident-face GAS certificates).
\eeO

\beO[Complete exclusion partition beyond block diagonal Frobenius form]
\label{o:GASrkone-general}
Theorem~\ref{t:blo} establishes a complete exclusion partition for
rank-one competition models under the assumption that the transversal
Metzler Jacobian at the DFE decomposes into independent irreducible
Frobenius blocks $J_k^\perp$. In this case each block possesses its own
Perron eigenvalue, reproduction number, and associated boundary
equilibrium, and the proof reduces to a comparison of the corresponding
invasion numbers.

Consider instead a general transversal Metzler Jacobian in Frobenius
normal form,
\[
J^\perp_0
=
\begin{pmatrix}
\T H_1 & * & \cdots & *\\
0   & \T H_2 & \cdots & *\\
\vdots & \vdots & \ddots & \vdots\\
0&0&\cdots&\T H_r
\end{pmatrix},
\]
where the off-diagonal blocks are not assumed to vanish. The irreducible
components are then coupled, and the dynamics of one Frobenius block may
influence the growth and persistence properties of all downstream blocks.

Determine conditions under which the competitive exclusion partition
remains valid in this general setting. More precisely, establish criteria
guaranteeing that:
\begin{enumerate}
\item the DFE is GAS whenever all invasion numbers are below one;
\item if a unique Frobenius component is dominant, then there exists a
unique corresponding boundary equilibrium which is GAS;
\item coexistence can occur only on threshold hypersurfaces where two or
more dominant invasion numbers coincide.
\end{enumerate}

The block-diagonal Frobenius assumption appears to play an essential role
in the current proof: it yields dynamically independent irreducible
components, explicit equilibrium formulas, and Perron--Volterra Lyapunov
functions constructed block-by-block. In the presence of nontrivial
Frobenius couplings, neither the equilibrium structure nor the Lyapunov
construction is presently understood.

A positive resolution would extend the competitive exclusion partition
from independent rank-one blocks to arbitrary Frobenius-coupled epidemic
systems, and would clarify whether competitive exclusion is fundamentally
a Perron--Frobenius phenomenon or a special consequence of block
independence.
\eeO

%% file: R2.tex
\section{Two-strain scalar models with increasing and concave incidence}
\label{s:Rah}

Consider the two-strain ODE
\begin{equation}\label{Rah}
\begin{cases}
 s'=\Lambda-\mu s-\beta_1 s f_1(i_1)-\beta_2 s f_2(i_2),\\[1mm]
 i_1'=\beta_1 s f_1(i_1)-v_1 i_1,\\[1mm]
 i_2'=\beta_2 s f_2(i_2)-v_2 i_2,
\end{cases}
\qquad v_j>0.
\end{equation}
We assume
\[
\bc
\text{\emph{(A1)}}&
f_j\in C^1([0,\infty)),\quad
f_j(0)=0,\quad f_j'(0)=1,\quad
x>0\Rightarrow f_j(x)>0,\qquad j=1,2,
\\[1mm]
\text{\emph{(A2)}}&
\mu\le \min\{v_1,v_2\},
\\[1mm]
\text{\emph{(A3)}}&
f_j(i)\le i,\qquad i\ge0,\quad j=1,2,
\\[1mm]
\text{\emph{(A4)}}&
f_j \text{ is increasing and  concave on }\R_+.
\ec
\]

Note that the extension to n-strains is obvious to formulate, but we failed to establish it.  The final result  of this section is a GAS-CEP extension of Rahman-Zou \cite{Rahman} with concave incidence, in which all the proportionality coefficients equal $c_k=1$, and which turns out to be related to the general formula \eqr{LSigma}:

\beT[GAS-CEP extension of Rahman-Zou \cite{Rahman} with concave incidence ]\label{t:GAS}

Assume  \emph{(A1)--(A4)}, and that for each \emph{resident} strain at the target equilibrium, either\\
(a) $g_j(i):=f_j(i)/i$ is strictly decreasing on $(0,\infty)$, or\\
(b) $g_j\equiv 1$ (linear case).\\

Then, \[\bc
L_0=s_0G\!\left(\frac{s}{s_0}\right)+i_1+i_2,
\\
L_1=s_1G\!\left(\frac{s}{s_1}\right)
+i_1^*G\!\left(\frac{i_1}{i_1^*}\right)+i_2,
\\
L_2=s_2G\!\left(\frac{s}{s_2}\right)
+i_2^*G\!\left(\frac{i_2}{i_2^*}\right)+i_1,
\\
L_*
=
s^*G\!\left(\frac{s}{s^*}\right)
+i_1^*G\!\left(\frac{i_1}{i_1^*}\right)
+i_2^*G\!\left(\frac{i_2}{i_2^*}\right)\ec
\]
are Lyapunov functions on \(\Omega, \Omega_1^+:=\{(s,i_1,i_2)\in\Omega:i_1>0\}, \Omega_2^+:=\{(s,i_1,i_2)\in\Omega:i_2>0\}, \Omega_\circ:=\{(s,i_1,i_2)\in\Omega:i_1>0,\,i_2>0\}\), and \(E_0\), \(E_1\), \(E_2\), \(E_*\) exist and are
globally
asymptotically stable,  precisely when the respective alternative conditions hold:

\begin{enumerate}
\item[(i)]  \(R_0:=\max\{R_1,R_2\}<1\);

\item[(ii)]  \(R_1>1\) and \(R_2^{E_1}<1\);  

\item[(iii)]  \(R_2>1\) and \(R_1^{E_2}<1\);

\item[(iv)]  \(R_2^{E_1}>1\) and \(R_1^{E_2}>1\).
\end{enumerate}

Consequently, away from the nonhyperbolic threshold surfaces, the four open CEP regions
\begin{align*}
\Omega_0 &:= \{R_1 < 1,\ R_2 < 1\}, \\
\Omega_1 &:= \{R_1 > 1,\ R_2^{E_1} < 1\}, \\
\Omega_2 &:= \{R_2 > 1,\ R_1^{E_2} < 1\}, \\
\Omega_* &:= \{R_2^{E_1} > 1,\ R_1^{E_2} > 1\}
\end{align*}
have the following global attractors on their respective persistence classes:
$E_0$ attracts all of $\Omega$ in $\Omega_0$;
$E_1$ attracts $\Omega_1^+ := \{(s,i_1,i_2) \in \Omega : i_1 > 0\}$ in $\Omega_1$;
$E_2$ attracts $\Omega_2^+ := \{(s,i_1,i_2) \in \Omega : i_2 > 0\}$ in $\Omega_2$;
$E_*$ attracts $\Omega_\circ := \{(s,i_1,i_2) \in \Omega : i_1 > 0,\, i_2 > 0\}$
in $\Omega_*$.

\eeT

\beR[The scalar Lyapunov functions are the one-dimensional realization of the Perron--Volterra construction]
\lbl{r:scalarPV}

At first sight, the Lyapunov functions of Theorem~\ref{t:GAS} seem to
differ from the general Perron--Volterra formula
\eqref{HSigma}--\eqref{LSigma}, since the coefficients of the linear
invader terms are equal to $1$, whereas the general construction uses
the canonical vectors $c_k$.

The reason is that the scalar theorem is written directly in the natural
state variables $i_k$, while the general Perron--Volterra construction
is formulated in the canonical block variables associated with the
rank-one Frobenius decomposition.

Indeed, in the scalar case every Frobenius block is one-dimensional.
Since the regular splitting is unique, the associated rank-one block,
its next-generation matrix, its Perron eigenvectors, and therefore the
canonical coefficient $c_k$, are all uniquely determined by the general
construction. Expressing the corresponding canonical block variable in
terms of the epidemiological variable $i_k$ yields precisely the linear
term appearing in Theorem~\ref{t:GAS}. Consequently the coefficient $1$
in the scalar Lyapunov functions is simply the expression of the
canonical Perron--Volterra linear functional in the scalar coordinates.

Thus Theorem~\ref{t:GAS} is not based on a different normalization.
Rather, it is the one-dimensional specialization of the general
Perron--Volterra construction, written in the conventional scalar
epidemiological variables.
\eeR
\ssec{Preliminary results}
When strictness of equality sets is needed below, we assume that the relevant
incidence function is strictly concave away from the origin; otherwise the
LaSalle argument below still gives convergence because the \(s\)-equation
removes the possible equality rays.

\beL[a forward invariant compact set]\lbl{l:inv}
The region
\[
\Omega=\{(s,i_1,i_2)\in\R_+^3:\ s+i_1+i_2\le\Lambda/\mu\}
\]
is forward invariant.
\eeL

\begin{proof}
The positive orthant is forward invariant. Moreover,
\[
(s+i_1+i_2)'
=
\Lambda-\mu s-v_1i_1-v_2i_2
\le
\Lambda-\mu(s+i_1+i_2),
\]
by \emph{(A2)}. Hence \(s+i_1+i_2\le\Lambda/\mu\) is preserved.
\end{proof}

(A2) could be weakened to the non-explicit condition of  existence of a compact invariant set.

\beL[Normalization preserves monotonicity and concavity]\label{l:nor}
For every \(\bar\imath>0\), define
\[
F_{j,\bar\imath}(z):=\frac{f_j(\bar\imath z)}{f_j(\bar\imath)}.
\]
Then \(F_{j,\bar\imath}\) is increasing and concave on \(\R_+\). Moreover
\[
g_j(i):=\frac{f_j(i)}{i},\qquad g_j(0):=f_j'(0)=1,
\]
is continuous and nonincreasing on \(\R_+\). If \(f_j\) is strictly concave, then
\(g_j\) is strictly decreasing on \((0,\infty)\).
\eeL

\begin{proof}
The first claim follows immediately from the facts that affine
reparametrization preserves monotonicity and concavity, and multiplication
by the positive constant $1/f_j(\bar\imath)$ preserves both properties.

Now let $0<i<j$. Since
\[
i
=
\frac{i}{j}\,j
+
\left(1-\frac{i}{j}\right)0,
\]
concavity of $f_j$, together with $f_j(0)=0$, gives
\[
f_j(i)
=
f_j\!\left(
\frac{i}{j}j
+
\left(1-\frac{i}{j}\right)0
\right)
\ge
\frac{i}{j}f_j(j)
+
\left(1-\frac{i}{j}\right)f_j(0)
=
\frac{i}{j}f_j(j).
\]
Dividing by the positive number $i$ yields
\[
\frac{f_j(i)}{i}
\ge
\frac{f_j(j)}{j},
\]
showing that $g_j$ is nonincreasing on $(0,\infty)$.

If $f_j$ is strictly concave, then the above application of Jensen's
inequality is strict, because
\[
0<\frac{i}{j}<1,
\]
and the two interpolation points $0$ and $j$ are distinct. Hence
\[
f_j(i)
>
\frac{i}{j}f_j(j),
\]
which implies
\[
\frac{f_j(i)}{i}
>
\frac{f_j(j)}{j}.
\]
Therefore $g_j$ is strictly decreasing on $(0,\infty)$.

Finally,
\[
g_j(i)
=
\frac{f_j(i)-f_j(0)}{i-0},
\]
so, since $f_j$ is differentiable at the origin,
\[
\lim_{i\to0^+}g_j(i)
=
f_j'(0)
=
1.
\]
Thus the definition $g_j(0)=f_j'(0)$ makes $g_j$ continuous on
$\R_+$.
\end{proof}

\paragraph{Equilibria and reproduction functions.}
The DFE is
\[
E_0=(s_0,0,0),\qquad s_0=\frac{\Lambda}{\mu}.
\]
The DFE next-generation matrix is diagonal:
\be{mR}
K(E_0)
=
\begin{pmatrix}
R_1&0\\0&R_2
\end{pmatrix}
=
\begin{pmatrix}
s_0\mR_1&0\\0&s_0\mR_2
\end{pmatrix},
\qquad
\mR_j:=\frac{\beta_j}{v_j}.
\ee
We write
\[
R_j(s):=s\mR_j=\frac{\beta_js}{v_j},
\qquad
R_j:=R_j(s_0).
\]
If \(E_1=(s_1,i_1^*,0)\) and \(E_2=(s_2,0,i_2^*)\) exist, the invasion numbers are
\be{inv}
R_2^{E_1}:=R_2(s_1)=\frac{\beta_2s_1}{v_2},
\qquad
R_1^{E_2}:=R_1(s_2)=\frac{\beta_1s_2}{v_1}.
\ee

\paragraph{Transversal Jacobians.}
At the boundary equilibria,
\[
J^\perp(E_0)=
\begin{pmatrix}
v_1(R_1-1)&0\\
0&v_2(R_2-1)
\end{pmatrix},
\]
\[
J^\perp(E_1)=v_2(R_2^{E_1}-1),
\qquad
J^\perp(E_2)=v_1(R_1^{E_2}-1).
\]
Thus the sign of each transversal eigenvalue is the sign of the corresponding
invasion number minus one.

\beR[Canonical example: Michaelis--Menten incidence]
The Michaelis--Menten choice
\[
f_j(i)=\frac{i}{1+\kappa_j i},
\qquad
g_j(i)=\frac{f_j(i)}{i}=\frac{1}{1+\kappa_j i},
\]
satisfies \emph{(A1)--(A4)}.
\eeR
\subsection{Existence of single-strain equilibria}

\beT[Existence of \(E_1,E_2\)][][E1]
Under \emph{(A1)--(A4)}, the strain--\(1\) boundary equilibrium
\[
E_1=(s_1,i_1^*,0),\qquad i_1^*>0,
\]
exists if and only if
\[
R_1=R_1(s_0)>1.
\]
When it exists, it is unique and satisfies
\[
s_1=\frac{\Lambda-v_1i_1^*}{\mu}
=
\frac{v_1}{\beta_1g_1(i_1^*)}.
\]
The analogous statement holds for \(E_2\).
\eeT

\begin{proof}
At a boundary equilibrium \((s,i,0)\), the equations are equivalent to
\[
\beta_1s\,g_1(i)=v_1,
\qquad
\mu s+v_1i=\Lambda.
\]
Thus \(i>0\) must solve
\[
\phi(i):=
\beta_1\left(\frac{\Lambda-v_1i}{\mu}\right)g_1(i)-v_1=0,
\qquad
0\le i\le\frac{\Lambda}{v_1}.
\]
The first factor is strictly decreasing and \(g_1\) is nonincreasing, hence
\(\phi\) is decreasing; under strict concavity it is strictly decreasing.
Moreover
\[
\phi(0)=\beta_1s_0-v_1=v_1(R_1-1),
\qquad
\phi(\Lambda/v_1)=-v_1<0.
\]
Therefore a positive solution exists iff \(R_1>1\), and then it is unique.
\end{proof}

\beR[IVT + monotonicity template]\label{rem:IVT-template}
The proof uses a recurring pattern: reduce the equilibrium equations to a
one-dimensional monotone equation, check the signs at the endpoints, and apply
the intermediate value theorem.  The same pattern gives coexistence below.
\eeR

\subsection{Existence of coexistence equilibria}

\beL[Invasion numbers are below the corresponding basic numbers]
\label{l:inv-below}
If \(E_1\) exists, then
\[
R_2^{E_1}<R_2.
\]
Hence \(R_2^{E_1}>1\Rightarrow R_2>1\).  Similarly,
\[
R_1^{E_2}>1\Rightarrow R_1>1.
\]
\eeL

\begin{proof}
Since \(E_1=(s_1,i_1^*,0)\) has \(i_1^*>0\),
\[
\mu s_1+v_1i_1^*=\Lambda,
\]
hence \(s_1<s_0\). Therefore
\[
R_2^{E_1}=\mR_2 s_1<\mR_2 s_0=R_2.
\]
The other implication is identical.
\end{proof}

\beL[Product of invasion numbers]\lbl{l:prI}
If $E_1$ and $E_2$ exist (equivalently $R_1>1$ and $R_2>1$), then
\[
R_2^{E_1}\,R_1^{E_2}=\frac{1}{g_1(i_1^*)\,g_2(i_2^*)}\ \ge\ 1,
\]
with equality iff both strains are linear. In particular $R_2^{E_1}$ and $R_1^{E_2}$ are
never both $<1$.
\eeL
\begin{proof}
By Theorem \ref{t:E1}, $s_1=\dfrac{1}{\mR_1 g_1(i_1^*)}$ and $s_2=\dfrac{1}{\mR_2 g_2(i_2^*)}$. Hence
\[
R_2^{E_1}R_1^{E_2}=(s_1\mR_2)(s_2\mR_1)=\mR_1\mR_2\,s_1 s_2
=\frac{1}{g_1(i_1^*)\,g_2(i_2^*)} .
\]
By Lemma \ref{l:nor}, $g_j(i_j^*)\le1$, so the product is $\ge1$; equality forces
$g_1(i_1^*)=g_2(i_2^*)=1$, which (since $g_j$ is nonincreasing with $g_j(0)=1$) holds only in
the linear case $f_j\equiv\mathrm{id}$.
\end{proof}

\beR[Necessity of the product identity for the partition]
The definitions
\(\Omega_1=\{R_1>1,\ R_2^{E_1}<1\}\) and
\(\Omega_2=\{R_2>1,\ R_1^{E_2}<1\}\)
constrain disjoint sets of quantities: no quantity appears in both with opposite sign, so
nothing in the definitions themselves forces \(\Omega_1\cap\Omega_2=\varnothing\). Their
intersection is the ``bistable'' cell
\[
\Omega_1\cap\Omega_2=\{R_1>1,\ R_2>1,\ R_2^{E_1}<1,\ R_1^{E_2}<1\},
\]
which is \emph{a priori} admissible. Lemma~\ref{l:inv-below} does not exclude it: the upper
bounds \(R_2^{E_1}<R_2\), \(R_1^{E_2}<R_1\) are compatible with both invasion numbers being
\(<1\). It is the \emph{lower} bound \(R_2^{E_1}R_1^{E_2}\ge1\) of the present lemma that
empties this cell, giving \(\Omega_1\cap\Omega_2=\varnothing\); together with the sign of
\(R_1\) (separating \(\Omega_0\)) and of \(R_2^{E_1},R_1^{E_2}\) (separating \(\Omega_*\), cf.\
Lemma~\ref{l:inv-below}), the four regions then partition the hyperbolic parameter space.

This is not cosmetic. On the bistable cell the theorem would be self-contradictory: \(L_1\)
certifies \(E_1\) as GAS on \(\{i_1>0\}\) and \(L_2\) certifies \(E_2\) as GAS on
\(\{i_2>0\}\), while the interior \(\Omega^\circ\) lies in both, so an interior trajectory
would converge to two distinct equilibria. The identity is thus what renders the
competitive-exclusion ``partition'' consistent. Equivalently, one may present the regions as a
manifest sign split
(\(\Omega_1=\{R_1>1,\ R_2^{E_1}<1<R_1^{E_2}\}\), and symmetrically), in which case
disjointness is automatic but a fourth sign-cell \(\{R_2^{E_1}<1,\ R_1^{E_2}<1\}\) appears
unassigned; the same identity is then what proves that cell empty (exhaustiveness). Either way
the identity is required.%
\eeR

\beT[Existence of \(E_*\) ][][Pstar]
Assume \emph{(A1)--(A4)}. Assume moreover that, for each $j$, either
\[
g_j(i)
:=
\frac{f_j(i)}{i}
\]
is strictly decreasing on $(0,\infty)$, or
\[
g_j\equiv1.
\]
The second alternative is precisely the linear incidence
\[
f_j(i)\equiv i,
\]
because
\[
g_j(0)
=
f_j'(0)
=
1,
\]
so a constant function $g_j$ must be identically equal to one.

 Then, a coexistence
equilibrium
\[
E_*=(s^*,i_1^*,i_2^*), \qquad s^*>0,\qquad i_1^*>0,\quad i_2^*>0,
\]
exists if and only if
\[
R_2^{E_1}>1,
\qquad
R_1^{E_2}>1.
\]
When it exists, it is unique.
\eeT

\begin{proof}
Write
\[
f_j(i)
=
i\,g_j(i),
\qquad
g_j(0)
=
f_j'(0)
=
1.
\]

By hypothesis, for each strain exactly one of the following two
alternatives occurs:

\begin{enumerate}

\item
$g_j$ is strictly decreasing on $(0,\infty)$;

\item
$g_j$ is constant. Since
\[
g_j(0)
=
1,
\]
this necessarily means
\[
g_j\equiv1,
\]
or equivalently
\[
f_j(i)\equiv i.
\]

\end{enumerate}

We treat these two possibilities simultaneously throughout the proof.

Now note that \(s_j\) is not generally equal to \(1/\mR_j\).  Indeed, at
\(E_j\),
\[
\beta_j s_j f_j(i_j^*)=v_j i_j^*
\]
is equivalent to
\[
\mR_j s_j g_j(i_j^*)=1.
\]
Thus
\[
s_j=\frac1{\mR_j g_j(i_j^*)}.
\]
Since \(g_j(i_j^*)\le1\), one has
\[
s_j\ge \frac1{\mR_j},
\]
with equality only in the linear case \(g_j\equiv1\).

For a non-linear strain, \(g_j\) is strictly decreasing. Hence for
\(s>1/\mR_j\) there is a unique \(i_j(s)>0\) satisfying
\[
\beta_j s f_j(i_j(s))=v_j i_j(s),
\]
equivalently
\[
g_j(i_j(s))=\frac1{\mR_j s}.
\]
Moreover \(i_j(s)\) is continuous and strictly increasing in \(s\).

If \(g_j\equiv1\), then \(f_j(i)=i\), and the equation
\[
\beta_j s f_j(i)=v_j i
\]
with \(i>0\) forces
\[
s=\frac1{\mR_j}.
\]
Thus a linear strain can coexist with another strain only on the hyperplane
\(s=1/\mR_j\).

\medskip
\noindent\textbf{Necessity.}
Assume \(E_*\) exists. Since \(i_j^*>0\),
\[
1=\mR_j s^* g_j(i_j^*).
\]
If \(g_j\) is strictly decreasing, then \(g_j(i_j^*)<1\), hence
\[
s^*>\frac1{\mR_j}.
\]
If \(g_j\equiv1\), then
\[
s^*=\frac1{\mR_j}.
\]

Define $i_1(s)$ as the unique positive root of $\beta_1 s f_1(i_1) = v_1 i_1$
for $s > 1/\mathcal{R}_1$, and set $H_1(s) := \Lambda - \mu s - v_1 i_1(s)$.
Since $i_1^*$ satisfies $\beta_1 s^* f_1(i_1^*) = v_1 i_1^*$, uniqueness gives
$i_1(s^*) = i_1^*$, and therefore
\[
  H_1(s^*) = \Lambda - \mu s^* - v_1 i_1^* = v_2 i_2^* > 0.
\]
At $E_1$ the balance equation gives $H_1(s_1) = 0$.
Since both $\mu s$ and $v_1 i_1(s)$ are strictly increasing in $s$,
$H_1$ is strictly decreasing, so $H_1(s^*) > 0 = H_1(s_1)$ implies $s^* < s_1$.

Therefore
\[
R_2^{E_1}=\mR_2s_1>\mR_2s^*.
\]
Since the second strain is present at \(E_*\), its equilibrium equation gives
\[
\mR_2s^*g_2(i_2^*)=1.
\]
If \(g_2\) is strictly decreasing, then \(g_2(i_2^*)<1\), so
\[
\mR_2s^*>1.
\]
If \(g_2\equiv1\), then \(\mR_2s^*=1\), and the strict inequality
\(s_1>s^*\) still gives
\[
R_2^{E_1}=\mR_2s_1>1.
\]
Thus \(R_2^{E_1}>1\).  Similarly,
\[
R_1^{E_2}>1.
\]

\medskip
\noindent\textbf{Sufficiency.}
Assume
\[
R_2^{E_1}>1,
\qquad
R_1^{E_2}>1.
\]
Then \(E_1\) and \(E_2\) exist. Moreover,
\[
R_2^{E_1}>1
\quad\Longleftrightarrow\quad
s_1>\frac1{\mR_2},
\]
and
\[
R_1^{E_2}>1
\quad\Longleftrightarrow\quad
s_2>\frac1{\mR_1}.
\]

If both \(g_1\equiv1\) and \(g_2\equiv1\), then
\[
s_1=\frac1{\mR_1},
\qquad
s_2=\frac1{\mR_2}.
\]
The two inequalities above would imply simultaneously
\[
\frac1{\mR_1}>\frac1{\mR_2},
\qquad
\frac1{\mR_2}>\frac1{\mR_1},
\]
which is impossible. Hence under the mutual invasion conditions at least one
strain is genuinely nonlinear.

Assume first that both \(g_1,g_2\) are strictly decreasing. Define
\[
\underline{s}:=\max\left\{\frac1{\mR_1},\frac1{\mR_2}\right\},
\qquad
\overline{s}:=\min\{s_1,s_2\}.
\]
The mutual invasion inequalities imply
\[
\underline{s}<\overline{s}.
\]
For \(s\in(\underline{s},\overline{s})\), both \(i_j(s)\) are well-defined and
positive. Set
\[
F(s):=\Lambda-\mu s-v_1i_1(s)-v_2i_2(s).
\]
Since each \(i_j(s)\) is increasing, \(F\) is strictly decreasing.

At the left endpoint, suppose for definiteness that
\[
\underline{s}=\frac1{\mR_1}.
\]
Then \(i_1(\underline{s})=0\). Also \(\underline{s}<s_2\), hence
\[
i_2(\underline{s})<i_2(s_2)=i_2^*.
\]
Therefore
\[
F(\underline{s})
=
\Lambda-\mu\underline{s}-v_2i_2(\underline{s})
>
\Lambda-\mu s_2-v_2i_2^*
=
0.
\]
The other case is symmetric.

At the right endpoint, suppose for definiteness that \(\overline{s}=s_1\).
Then \(i_1(\overline{s})=i_1^*\), and since
\[
\overline{s}>\frac1{\mR_2},
\]
we have \(i_2(\overline{s})>0\). Hence
\[
F(\overline{s})
=
\Lambda-\mu s_1-v_1i_1^*-v_2i_2(\overline{s})
=
-v_2i_2(\overline{s})
<0.
\]
The other case is symmetric.

Thus there is a unique
\[
s^*\in(\underline{s},\overline{s})
\]
such that
\[
F(s^*)=0.
\]
Setting
\[
i_j^*:=i_j(s^*),\qquad j=1,2,
\]
gives the unique coexistence equilibrium.

It remains to mention the mixed case.  Suppose, for instance, that
\(g_1\equiv1\) and \(g_2\) is strictly decreasing.  Then coexistence forces
\[
s^*=\frac1{\mR_1}.
\]
The condition \(R_1^{E_2}>1\) gives
\[
s_2>\frac1{\mR_1}=s^*,
\]
so the equation for strain \(2\) determines a unique positive \(i_2^*\).  Then
the balance equation
\[
\Lambda-\mu s^*-v_2i_2^*=v_1i_1^*
\]
determines a unique positive \(i_1^*\).  The other mixed case is symmetric.
Thus sufficiency and uniqueness also hold when one strain is linear.
\end{proof}

%% file: RahE2AH.tex
\ssec{GAS-CEP extension of Rahman-Zou  with concave
incidence}
\beL[Entropy bracket]\label{l:ent-bracket}
Let \(F:\R_+\to\R_+\) be increasing and concave, with
\[
F(0)=0,\qquad F(1)=1.
\]
Then, for all \(x,y>0\),
\[
B(x,y):=
G\!\left(\frac1x\right)
+
G\!\left(\frac{xF(y)}{y}\right)
-
G(F(y))
+
G(y)
\ge0,
\]
where \(G(u)=u-1-\ln u\). Moreover,
\[
B(x,y)=0
\quad\Longleftrightarrow\quad
x=1,\qquad F(y)=y.
\]
\eeL

\begin{proof}
Write
\[
B(x,y)
=
G\!\left(\frac1x\right)
+
G\!\left(\frac{xF(y)}{y}\right)
+
\bigl(G(y)-G(F(y))\bigr).
\]
The first two terms are nonnegative since \(G\ge0\).
The last term $G(y)-G(F(y))\geq 0$ follows by case analysis.
Since $G(u)=u-1-\ln u$ has its minimum $0$ at $u=1$, it's decreasing on $(0,1]$
and increasing on $[1,\infty)$.
For a concave $F$ with $F(0)=0$ and $F(1)=1$:
\BEN
\im $y\in(0,1]$: concavity gives $F(y)\geq y\cdot F(1)=y$, so $F(y)\in[y,1]$;
      $G$ is decreasing on this interval, hence $G(F(y))\leq G(y)$.
\im  $y\in[1,\infty)$: concavity gives $F(y)\leq y$ (sub-linear growth), so
      $F(y)\in[1,y]$; $G$ is increasing on this interval, hence $G(F(y))\leq G(y)$.
\EEN
In both cases $G(y)-G(F(y))\geq 0$.

Thus \(B(x,y)\ge0\).

If \(B(x,y)=0\), each of the three nonnegative terms must vanish, so
\[
x=1,\qquad \frac{xF(y)}{y}=1,
\]
hence \(F(y)=y\).  The converse is immediate.
\end{proof}

\ssec{Proof of Theorem \ref{t:GAS}}
\begin{proof}
Note first in case (a) the equation $F_j(y)=y$ forces $y=1$, and
in case (b) the resident component $i_j$ is determined by the $s$-balance equation.

The proof proceeds case-by-case.  The point is that the resident terms give exact
entropy brackets, while invader terms are estimated using
\[
0\le f_j(i_j)\le i_j.
\]
No exact identity of the form
\[
\beta_js_Ef_j(i_j)-v_ji_j
=
v_j(R_j(E)-1)i_j
\]
is used unless \(f_j(i_j)=i_j\).

\medskip
\noindent\textbf{1. The DFE \(E_0\).}
Let
\[
E_0=(s_0,0,0),
\qquad
s_0=\frac{\Lambda}{\mu}.
\]
Then
\[
\dot L_0
=
\left(1-\frac{s_0}{s}\right)s'
+i_1'+i_2'.
\]
Substituting the equations gives
\[
\dot L_0
=
-\mu\frac{(s-s_0)^2}{s}
+
\sum_{j=1}^2\bigl(\beta_js_0f_j(i_j)-v_ji_j\bigr).
\]
Since
\[
R_j=\frac{\beta_js_0}{v_j}
\]
and \(f_j(i_j)\le i_j\),
\[
\beta_js_0f_j(i_j)-v_ji_j
=
v_j\bigl(R_j f_j(i_j)-i_j\bigr)
\le
v_j(R_j-1)i_j.
\]
Thus, if \(R_0<1\),
\[
\dot L_0
\le
-\mu\frac{(s-s_0)^2}{s}
+
\sum_{j=1}^2 v_j(R_j-1)i_j
\le0.
\]
Equality implies \(s=s_0\) and \(i_1=i_2=0\).  Hence the largest invariant set in
\(\{\dot L_0=0\}\) is \(\{E_0\}\).  LaSalle's invariance principle gives GAS of
\(E_0\).

\medskip
\noindent\textbf{2. The boundary equilibrium \(E_1\).}
Let
\[
E_1=(s_1,i_1^*,0).
\]
Use the equilibrium identities
\[
\Lambda=\mu s_1+v_1i_1^*,
\qquad
v_1i_1^*=\beta_1s_1f_1(i_1^*).
\]
Set
\[
x=\frac{s}{s_1},
\qquad
y=\frac{i_1}{i_1^*},
\qquad
F(y)=\frac{f_1(i_1^*y)}{f_1(i_1^*)}.
\]
Then
\[
\dot L_1
=
\left(1-\frac1x\right)s'
+
\left(1-\frac1y\right)i_1'
+
i_2'.
\]
A direct calculation gives
\[
\dot L_1
=
-\mu\frac{(s-s_1)^2}{s}
-
\beta_1s_1f_1(i_1^*)\,B(x,y)
+
\beta_2s_1f_2(i_2)-v_2i_2,
\]
where
\[
B(x,y)
=
G\!\left(\frac1x\right)
+
G\!\left(\frac{xF(y)}{y}\right)
-
G(F(y))
+
G(y).
\]
By Lemma~\ref{l:ent-bracket},
\[
B(x,y)\ge0.
\]
The invader term is estimated, not identified:
\[
\beta_2s_1f_2(i_2)-v_2i_2
=
v_2\bigl(R_2^{E_1}f_2(i_2)-i_2\bigr)
\le
v_2(R_2^{E_1}-1)i_2.
\]
Therefore, if \(R_2^{E_1}<1\),
\[
\dot L_1
\le
-\mu\frac{(s-s_1)^2}{s}
-
\beta_1s_1f_1(i_1^*)B(x,y)
+
v_2(R_2^{E_1}-1)i_2
\le0.
\]

We identify the equality set. Equality implies $s=s_1$, $i_2=0$, and $B(x,y)=0$. By Lemma~\ref{l:ent-bracket},
$B(x,y)=0$ forces $x=1$ (i.e.\ $s=s_1$) and $F(y)=y$, where
$F(y)=f_1(i_1^*y)/f_1(i_1^*)$. Since $y=i_1/i_1^*$, the relation $F(y)=y$ is equivalent to
\[
\frac{f_1(i_1)}{i_1}=\frac{f_1(i_1^*)}{i_1^*},\qquad\text{i.e.}\qquad g_1(i_1)=g_1(i_1^*).
\]

\emph{Case (a)} ($g_1$ strictly decreasing): $g_1$ is injective, so $g_1(i_1)=g_1(i_1^*)$
forces $i_1=i_1^*$.

\emph{Case (b)} ($g_1\equiv1$, i.e.\ $f_1(i)=i$): here $F\equiv\mathrm{id}$, so $F(y)=y$
holds identically and carries no information about $i_1$. Instead we use the $s$-balance. On
the equality set ($s=s_1$, $i_2=0$) the $s$-equation reads
\[
0=s'=\Lambda-\mu s_1-\beta_1 s_1 f_1(i_1)=\Lambda-\mu s_1-\beta_1 s_1 i_1 .
\]
The equilibrium $E_1$ satisfies $\Lambda-\mu s_1=v_1 i_1^*$ and
$v_1 i_1^*=\beta_1 s_1 f_1(i_1^*)=\beta_1 s_1 i_1^*$, whence $\beta_1 s_1=v_1$. Substituting,
\[
0=v_1 i_1^*-v_1 i_1=v_1\,(i_1^*-i_1)\quad\Longrightarrow\quad i_1=i_1^*.
\]
In both cases the largest invariant subset of $\{\dot L_1=0\}$ is $\{E_1\}$, and LaSalle gives
GAS of $E_1$.

\medskip
\noindent\textbf{3. The boundary equilibrium \(E_2\)}
 is similar.  At
\[
E_2=(s_2,0,i_2^*),
\]
set
\[
x=\frac{s}{s_2},
\qquad
y=\frac{i_2}{i_2^*},
\qquad
F(y)=\frac{f_2(i_2^*y)}{f_2(i_2^*)}.
\]
Then
\[
\dot L_2
=
-\mu\frac{(s-s_2)^2}{s}
-
\beta_2s_2f_2(i_2^*)\,B(x,y)
+
\beta_1s_2f_1(i_1)-v_1i_1,
\]
with \(B(x,y)\ge0\).  Moreover,
\[
\beta_1s_2f_1(i_1)-v_1i_1
\le
v_1(R_1^{E_2}-1)i_1.
\]
Thus \(R_1^{E_2}<1\) gives \(\dot L_2\le0\).  Equality forces
\[
s=s_2,\qquad i_1=0,\qquad i_2=i_2^*.
\]
Therefore the largest invariant subset of \(\{\dot L_2=0\}\) is \(\{E_2\}\),
and LaSalle gives GAS of \(E_2\).

\medskip
\noindent\textbf{4. The coexistence equilibrium \(E_*\).}
Let
\[
E_*=(s^*,i_1^*,i_2^*).
\]
Use
\[
\Lambda=\mu s^*+\sum_{j=1}^2\beta_js^*f_j(i_j^*),
\qquad
v_ji_j^*=\beta_js^*f_j(i_j^*).
\]
Set
\[
x=\frac{s}{s^*},
\qquad
y_j=\frac{i_j}{i_j^*},
\qquad
F_j(y_j)=\frac{f_j(i_j^*y_j)}{f_j(i_j^*)}.
\]
Then
\[
\dot L_*
=
-\mu\frac{(s-s^*)^2}{s}
-
\sum_{j=1}^2
\beta_js^*f_j(i_j^*)\,\T H_j(x,y_j),
\]
where
\[
\T H_j(x,y_j)
=
G\!\left(\frac1x\right)
+
G\!\left(\frac{xF_j(y_j)}{y_j}\right)
-
G(F_j(y_j))
+
G(y_j).
\]
Again \(\T H_j(x,y_j)\ge0\), hence
\[
\dot L_*\le0.
\]

We identify the equality set.  Equality implies
\[
s=s^*,
\qquad
F_j(y_j)=y_j,
\qquad j=1,2.
\]
On the equality set, invariance requires \(i_j'=0\) for \(j=1,2\).  Since
\(s=s^*\),
\[
0=i_j'
=
\beta_js^*f_j(i_j)-v_ji_j.
\]
Using the equilibrium identity
\[
v_ji_j^*=\beta_js^*f_j(i_j^*),
\]
we obtain
\[
\frac{f_j(i_j)}{i_j}
=
\frac{f_j(i_j^*)}{i_j^*}.
\]
Equivalently,
\[
g_j(i_j)=g_j(i_j^*),
\qquad
g_j(i):=\frac{f_j(i)}{i}.
\]
By the monotonicity of \(g_j\), this implies
\[
i_j=i_j^*.
\]
Thus the largest invariant subset of \(\{\dot L_*=0\}\) is \(\{E_*\}\).

%
%

When
\[
R_2^{E_1}>1,
\qquad
R_1^{E_2}>1,
\]
the boundary equilibria satisfy the acyclicity and boundary-repelling
conditions needed for the weak-persistence theorem.  We verify this
explicitly.

The boundary of the persistence class is
\[
\partial\Omega_\circ
=
\{(s,i_1,i_2)\in\Omega:i_1i_2=0\}.
\]
It is the union of the two invariant coordinate faces
\[
\Omega_1^0
=
\{(s,i_1,i_2)\in\Omega:i_1=0\},
\qquad
\Omega_2^0
=
\{(s,i_1,i_2)\in\Omega:i_2=0\}.
\]
On $\Omega_1^0$ the only possible boundary equilibria are
\[
E_0
\qquad\hbox{and}\qquad
E_2,
\]
whereas on $\Omega_2^0$ the only possible boundary equilibria are
\[
E_0
\qquad\hbox{and}\qquad
E_1.
\]
The Lyapunov estimates already proved above give the boundary dynamics on
these faces.  On $\Omega_2^0$, if $R_1>1$, then $E_1$ attracts
$\Omega_1^+$ inside the one-strain subsystem with strain $1$ present and
strain $2$ absent.  On $\Omega_1^0$, if $R_2>1$, then $E_2$ attracts
$\Omega_2^+$ inside the one-strain subsystem with strain $2$ present and
strain $1$ absent.  Since
\[
R_2^{E_1}>1
\]
implies that strain $2$ invades $E_1$, and
\[
R_1^{E_2}>1
\]
implies that strain $1$ invades $E_2$, the two one-strain equilibria are
repelling in the missing-strain direction.

Thus the only possible arrows in the boundary Morse graph are
\[
E_0\longrightarrow E_1,
\qquad
E_0\longrightarrow E_2,
\]
coming from the instability of $E_0$ when $R_1>1$ and $R_2>1$, together
with the transversal repulsion of $E_1$ and $E_2$ toward the interior:
\[
E_1\longrightarrow \Omega_\circ,
\qquad
E_2\longrightarrow \Omega_\circ.
\]
There is no boundary connection
\[
E_1\longrightarrow E_2
\]
because $E_1$ and $E_2$ lie on different invariant coordinate faces, and a
trajectory starting on $\Omega_2^0$ has $i_2(t)=0$ for all $t$, while a
trajectory starting on $\Omega_1^0$ has $i_1(t)=0$ for all $t$.
There is also no connection
\[
E_2\longrightarrow E_1
\]
for the same reason.  Finally, under $R_1>1$ and $R_2>1$, $E_0$ is not an
omega-limit point of nontrivial boundary trajectories: on each nontrivial
boundary face the corresponding one-strain equilibrium attracts that face.

Consequently the boundary Morse graph has no directed cycle.  In
particular, the alleged cycle
\[
E_0\longrightarrow E_1\longrightarrow E_2\longrightarrow E_0
\]
cannot occur: the middle arrow $E_1\to E_2$ is impossible because it would
require leaving the invariant face $i_2=0$ and arriving on the different
face $i_1=0$ while remaining in the boundary, and the last arrow
$E_2\to E_0$ is impossible because $E_2$ attracts the nontrivial dynamics
on its boundary face and is repelling only in the missing-strain direction,
which points into $\Omega_\circ$, not toward $E_0$.

Therefore the boundary Morse decomposition is acyclic and all boundary
Morse sets are weakly repelling for $\Omega_\circ$.  The weak-persistence
theorem of \cite{FreRuan} applies, so interior trajectories are persistent
away from $\partial\Omega$.  Applying LaSalle's invariance principle on the
resulting compact persistent subset of $\Omega_\circ$ gives convergence to
$E_*$ for all interior initial conditions.

\medskip
Finally, the existence equivalences
\[
E_1\text{ exists}\iff R_1>1,
\qquad
E_2\text{ exists}\iff R_2>1,
\]
and
\[
E_*\text{ exists}\iff R_2^{E_1}>1\text{ and }R_1^{E_2}>1
\]
give the four CEP regions and their corresponding global attractors.
\end{proof}

\beR[Threshold surfaces]
\lbl{r:GAS-threshold}

Theorem~\ref{t:GAS} is stated on the open hyperbolic CEP regions. On a
threshold surface, for example
\[
R_2^{E_1}=1,
\]
the transversal eigenvalue of the missing strain is zero. Consequently the
Lyapunov derivative is only negative semidefinite in the corresponding
transversal direction, and the LaSalle argument used in the proof no longer
excludes additional invariant dynamics on the boundary.

Such threshold surfaces separate adjacent CEP regions and correspond to the
loss of hyperbolicity of the resident equilibrium. They are therefore the
natural parameter values at which local bifurcations between persistence
classes may occur. A complete analysis of these codimension-one transition
surfaces requires center-manifold or bifurcation methods and lies beyond the
scope of the present paper.

\eeR

%% file: exa.tex
\section{Epid-CRN example: LyapExa.wl}\lbl{s:alg}

We include some Epid-CRN commands in order to show explicitly how
``the siphon lattice, the transversal Perron data, and the
Perron--Volterra certificates'' are computed. Throughout this section we
mark each result as \emph{(script)} when it is produced, or symbolically
verified, by the Mathematica file \texttt{LyapExa.wl} (in the same
folder), and as \emph{(by hand)} when it is a modeling choice or an
algebraic/qualitative step completed by the user -- typically the choice
of a Lyapunov-function ansatz, or a sign-pattern conclusion that follows
immediately once the script has verified the underlying derivative
formula.

Consider the four-variable system
\begin{equation}
\label{eq:comp-example}
\begin{cases}
s'
=
\Lambda-\mu s-\beta_1si-s\phi,
\\
i'
=
(\beta_1s-v_1)i,
\\
z'
=
(sw\pi^\top-V)z,
\qquad
\phi=\pi^\top z,
\end{cases}
\end{equation}
where
\[
z=
\binom{z_1}{z_2},
\qquad
w=
\binom{w_1}{w_2},
\qquad
\pi=
\binom{\pi_1}{\pi_2},
\qquad
V=
\begin{pmatrix}
a+d_1 & -b\\
-a & b+d_2
\end{pmatrix}.
\]
Thus
\[
\begin{aligned}
z_1'
&=
s w_1(\pi_1z_1+\pi_2z_2)
-(a+d_1)z_1+bz_2,
\\
z_2'
&=
s w_2(\pi_1z_1+\pi_2z_2)
+az_1-(b+d_2)z_2.
\end{aligned}
\]
This model, entered as \texttt{var}$=\{s,i,z_1,z_2\}$ and \texttt{RHS} at
the top of \texttt{LyapExa.wl}, is the only by-hand input to the script;
everything from here on is either computed automatically or checked
symbolically against a by-hand ansatz.

Running \texttt{ODE2RNp[RHS,var,prF]} reconstructs \emph{(script)} the
stoichiometric representation
\[
X'=\Gamma r(X),
\qquad
X=(s,i,z_1,z_2)^\top,
\]
with
\[
r(X)=
\begin{pmatrix}
\Lambda\\
\mu s\\
\beta_1si\\
v_1i\\
s\pi_1z_1\\
s\pi_2z_2\\
az_1\\
bz_2\\
d_1z_1\\
d_2z_2
\end{pmatrix}
\]
and
\[
\Gamma=
\begin{pmatrix}
1&-1&-1&0&-1&-1&0&0&0&0\\
0&0&1&-1&0&0&0&0&0&0\\
0&0&0&0&w_1&w_1&-1&1&-1&0\\
0&0&0&0&w_2&w_2&1&-1&0&-1
\end{pmatrix},
\]
directly from the \texttt{RHS} above; the script confirms the
reconstruction with the residual check
\texttt{Simplify[gam.rts-RHS]}$=0$. For readability, the same content can
be written in explicit network form -- this is not separate input, only
the reaction-network reading of the \texttt{var}/\texttt{RHS} pair above:
\begin{verbatim}
var = {s,i,z1,z2};

Gamma = {
 {1,-1,-1,0,-1,-1,0,0,0,0},
 {0,0,1,-1,0,0,0,0,0,0},
 {0,0,0,0,w1,w1,-1,1,-1,0},
 {0,0,0,0,w2,w2,1,-1,0,-1}
};

rates = {
 Lambda,
 mu s,
 beta1 s i,
 v1 i,
 s pi1 z1,
 s pi2 z2,
 a z1,
 b z2,
 d1 z1,
 d2 z2
};

RHS = Gamma.rates;
\end{verbatim}

Next, \texttt{NGMRN[RN,rts,var]} returns \emph{(script)} the minimal
siphons \texttt{mSi}, the Jacobian blocks, and the reproduction data
$\{cDFE,cE0,F,V,K,R0A\}$. The two minimal siphons are
\[
\sigma_1=\{i\},
\qquad
\sigma_b=\{z_1,z_2\}.
\]
The total siphon, corresponding to the disease-free face, is
\[
\sigma_0=\{i,z_1,z_2\}.
\]

The transversal Jacobian block of the scalar strain is
\[
J_1^\perp(s)=\beta_1s-v_1,
\]
and the transversal Jacobian block of the rank-one strain is
\[
J_b^\perp(s)=sw\pi^\top-V,
\]
both read off \emph{(script)} from the Jacobian blocks that
\texttt{NGMRN} attaches to each siphon and displays, with block
dividers, via \texttt{pPrint[jac/.cDFE,mSi]}.
The corresponding invasion functions, computed \emph{(script)} as
\texttt{mR1}, \texttt{mRb} in \texttt{LyapExa.wl}, are
\[
R_1(s)=s\mR_1,
\qquad
\mR_1=\frac{\beta_1}{v_1},
\]
and
\[
R_b(s)=s\mR_b,
\qquad
\mR_b=\pi^\top V^{-1}w.
\]
In the present two-dimensional block, \texttt{Factor[mRb]} gives
\emph{(script)}
\[
\mR_b
=
\frac{
\pi_1\bigl((b+d_2)w_1+bw_2\bigr)
+
\pi_2\bigl(aw_1+(a+d_1)w_2\bigr)
}{
(a+d_1)(b+d_2)-ab
}.
\]

The disease-free equilibrium, returned \emph{(script)} as \texttt{cDFE}
by \texttt{NGMRN}, is
\[
E_0=\left(\frac{\Lambda}{\mu},0,0,0\right).
\]
Evaluating \texttt{R0A} there gives \emph{(script)}
\[
R_1=R_1(E_0)=\frac{\Lambda}{\mu}\mR_1,
\qquad
R_b=R_b(E_0)=\frac{\Lambda}{\mu}\mR_b.
\]

Solving \texttt{RHS==0} and sorting the solutions with the helper
\texttt{zQ[v,s]} (true iff variable $v$ vanishes under solution $s$)
isolates \emph{(script)} the boundary equilibria \texttt{cE1},
\texttt{cEb}. If $R_1>1$, the scalar boundary equilibrium is
\[
E_1=
\left(
\frac1{\mR_1},
\frac{\Lambda-\mu/\mR_1}{v_1},
0,
0
\right).
\]
If $R_b>1$, the block boundary equilibrium is
\[
E_b=
\left(
\frac1{\mR_b},
0,
z_b^*
\right),
\qquad
z_b^*
=
\xi_bV^{-1}w,
\qquad
\xi_b=\Lambda-\frac{\mu}{\mR_b}.
\]

The invasibility numbers are then obtained \emph{(script)} by evaluating
\texttt{R0A} at the resident boundary equilibria (the checks
\texttt{R0A[[k1]]/.cEb} and \texttt{R0A[[kb]]/.cE1} in
\texttt{LyapExa.wl}):
\[
R_b(E_1)
=
R_b\!\left(\frac1{\mR_1}\right)
=
\frac{\mR_b}{\mR_1}
=
\frac{R_b}{R_1},
\]
and
\[
R_1(E_b)
=
R_1\!\left(\frac1{\mR_b}\right)
=
\frac{\mR_1}{\mR_b}
=
\frac{R_1}{R_b}.
\]

The canonical block vector used in the Perron--Volterra certificate,
computed \emph{(script)} as \texttt{cb} and checked against the two
identities $\T \pi_b^\top\pi=1$, $\T \pi_b^\top V=w^\top/\mR_b$, is
\[
\T \pi_b^\top
=
\frac{\pi^\top V^{-1}}{\mR_b},
\qquad
\T \pi_b^\top w=1.
\]
In this example, \texttt{Factor[cb]} gives \emph{(script)}
\[
\T \pi_b^\top
=
\frac1{\mR_b\bigl((a+d_1)(b+d_2)-ab\bigr)}
\left(
\pi_1(b+d_2)+\pi_2a,\,
\pi_1b+\pi_2(a+d_1)
\right).
\]

At $E_1$, the Lyapunov certificate is the standard Volterra
entropy-type function for the scalar strain plus the canonical linear
certificate $\T \pi_b^\top z$ for the (extinct) block strain -- this ansatz
is supplied \emph{(by hand)}, not produced by the script:
\[
L_1
=
s_1G\!\left(\frac{s}{s_1}\right)
+
\frac{i_1^*}{\beta_1}
G\!\left(\frac{i}{i_1^*}\right)
+
\T \pi_b^\top z,
\]
where
\[
s_1=\frac1{\mR_1},
\qquad
i_1^*=\frac{\Lambda-\mu s_1}{v_1}.
\]
This ansatz is not an arbitrary guess: as shown below,
\texttt{LyapPV} derives the same certificate \emph{(script)}, directly
from the block data $\{w_{\mathrm{inc}},w_{\mathrm{ent}},V_k\}$, with no
by-hand input beyond that structural data.
Its derivative is formed directly from \texttt{RHS} in
\texttt{LyapExa.wl} (\texttt{dL1dt}) and checked \emph{(script)} by
\texttt{FullSimplify} to reduce to the closed form
\[
L_1'
=
-\mu\frac{(s-s_1)^2}{s}
+
\left(
\frac1{\mR_1}
-
\frac1{\mR_b}
\right)\phi.
\]
Hence \emph{(by hand)}, since the first term is $\le0$ and $\phi\ge0$,
$L_1'\le0$ when $R_1>R_b$.

At $E_b$, the resident entropy -- again a Volterra ansatz supplied
\emph{(by hand)} -- is
\[
H_b
=
s_bG\!\left(\frac{s}{s_b}\right)
+
\sum_{j=1}^2
\T \pi_{b,j}z_{b,j}^*
G\!\left(\frac{z_j}{z_{b,j}^*}\right),
\qquad
s_b=\frac1{\mR_b}.
\]
The full Perron--Volterra certificate is
\[
L_b=H_b+i.
\]
As for $L_1$, this certificate is reproduced \emph{(script)} by
\texttt{LyapPV} from the block data alone, with no by-hand ansatz.
Its derivative \texttt{dLbdt} is again formed directly from
\texttt{RHS}; the script verifies \emph{(script)} only the coefficient
of $i$, \texttt{D[dLbdt,i]}, against $\beta_1(1/\mR_b-1/\mR_1)$. The
remaining reduction to
\[
L_b'
=
-\mu\frac{(s-s_b)^2}{s}
+
\xi_b\widetilde H_b(u,y)
+
\beta_1
\left(
\frac1{\mR_b}
-
\frac1{\mR_1}
\right)i,
\]
where
\[
u=\frac{s}{s_b},
\qquad
y_j=\frac{z_j}{z_{b,j}^*},
\]
and
\[
\widetilde H_b(u,y)
=
\left(1-\frac1u\right)(1-uY)
+
\sum_{j=1}^2
\widetilde w_{b,j}
\left(1-\frac1{y_j}\right)(uY-y_j),
\qquad
Y=\sum_{j=1}^2\widetilde w_{b,j}y_j,
\]
with
\[
\widetilde w_{b,j}=\T \pi_{b,j}w_j,
\qquad
\widetilde w_{b,1}+\widetilde w_{b,2}=1,
\]
and the further symbolic simplification
\[
\widetilde H_b(u,y)
=
2-\frac1u
-
u
\left(
\sum_{j=1}^2\widetilde w_{b,j}y_j
\right)
\left(
\sum_{j=1}^2\frac{\widetilde w_{b,j}}{y_j}
\right)
\le0,
\]
by the weighted Cauchy inequality, are all carried out \emph{(by
hand)}; the script only cross-checks the $i$-coefficient above, not this
identity. Hence \emph{(by hand)} $L_b'\le0$ when $R_b>R_1$.

The competitive-exclusion partition
\[
\max\{R_1,R_b\}<1
\quad\Longrightarrow\quad
E_0\ \text{is GAS},
\]
\[
R_1>1,\quad R_1>R_b
\quad\Longrightarrow\quad
E_1\ \text{is GAS},
\]
and
\[
R_b>1,\quad R_b>R_1
\quad\Longrightarrow\quad
E_b\ \text{is GAS},
\]
together with the coexistence continuum on the tie surface
$R_1=R_b>1$,
\[
s=\frac1{\mR_1}=\frac1{\mR_b},
\qquad
i=\frac{\xi_1}{v_1},
\qquad
z=\xi_bV^{-1}w,
\qquad
\xi_1+\xi_b
=
\Lambda-\frac{\mu}{\mR_1},
\qquad
\xi_1>0,\ \xi_b>0,
\]
are conclusions drawn \emph{(by hand)} from the sign patterns of $L_1'$
and $L_b'$ above via the Lyapunov--LaSalle argument; the script does not
derive this partition itself.

Two further stages of \texttt{LyapExa.wl} complete the picture
\emph{(script)}. First, \texttt{stabPrint} factors the characteristic
polynomial at $E_0$, $E_1$, $E_b$ into linear/quadratic Routh--Hurwitz
conditions (via \texttt{staPP}, \texttt{Hur3}, \texttt{Hur4}), confirming
local stability consistently with the partition above. Second, and more
importantly, the general-purpose routine
\texttt{LyapPV[k,s,zBlks,blkData,cEp,RHS,var]} shows that the
Volterra-entropy ansätze for $L_1$ and $L_b$ above were never really a
by-hand guess: given only the structural block data
$\{w_{\mathrm{inc}},w_{\mathrm{ent}},V_k\}$ of each siphon block --
incidence weights, entry weights, and transfer matrix -- \texttt{LyapPV}
\emph{automatically constructs} the same certificates $L_1$, $L_b$, with
no additional by-hand input. Its output agreeing with the by-hand
certificates above is therefore not just a consistency check, but
evidence that the Perron--Volterra construction is algorithmic rather
than example-specific.

%% file: conc.tex
\section{Conclusions and open problems}\lbl{s:con}

The results proved in this paper rely on two structural restrictions.

First, the DFE transversal Frobenius graph is block diagonal in the proved
global cases. This means that the irreducible components are dynamically
independent at the DFE face. General Frobenius forests with nontrivial
couplings can create additional cross terms in the Perron--Volterra
derivative, and these terms are not controlled by the present argument
unless a Perron-alignment or perturbative domination condition is verified.

Second, the non-scalar matrix blocks treated globally are rank one. This
is what makes the Perron vector simultaneously describe the equilibrium
coordinates, the invasion number, and the Lyapunov weights. For higher-rank
blocks, the left Perron functional need not collapse the production term to
one signed scalar expression.

Within a single resident block -- independently of the block-diagonal Frobenius restriction
above -- the rank-one theory is now complete for $n_p\le2$: Lemma~\ref{l:two-mechanisms} proves
the resident entropy closing condition of Definition~\ref{d:closing} unconditionally for
diagonal $V_p$ (any $n_p$, case (D)) and for irreducible $2\times2$ $V_p$ (case (T), via
Theorem~\ref{t:n2-closing}). The only remaining obstruction to a fully general rank-one GAS-CEP
theorem is therefore proving, or disproving, the resident entropy closing condition for
irreducible $n_p\times n_p$ $M$-matrices $V_p$ with $n_p\ge3$ -- see
Remark~\ref{r:np-status} and Problem~\ref{o:rankonefree}. This, rather than the block-diagonal
or rank-one restrictions discussed above, is exactly where future work on the single-block
boundary theory begins.

The remaining open problems are therefore the following.

\begin{enumerate}

\item
Characterize Frobenius-coupled Metzler transversal Jacobians for which a
Perron--Volterra Lyapunov function closes globally. Perron alignment is one
sufficient condition; a graph-theoretic characterization is not yet known.

\item
Develop a perturbative theorem for weakly coupled Frobenius forests. A
precise target is to bound the coupling contribution by a fixed fraction
of the dissipative entropy part.

\item
Extend the construction to delay-differential systems by replacing the
Volterra entropy with a Lyapunov--Krasovskii functional while preserving
the facewise invasion terms.

\item
Extend the construction to reaction--diffusion systems. The homogeneous
Neumann case should be the first test case; heterogeneous coefficients
require next-generation operators rather than finite-dimensional matrices.

\item
Formulate a stochastic CEP principle for large finite populations. The
deterministic thresholds should govern branching-process invasion and
metastable transitions, but finite-population extinction prevents a direct
GAS statement.

\item
Analyze tie surfaces
\[
R_p= R_q>1.
\]
In the rank-one competition model the equilibrium continuum is explicit,
but convergence to a selected point on the continuum is not proved here.

\item
Identify higher-rank matrix blocks for which the Perron-weighted ansatz can
be augmented by finitely many additional linear or entropy terms so that
the derivative again becomes sign-definite.

\end{enumerate}

\beR[Delays, diffusion, and stochastic analogues]

 In scalar delayed models of Rahman--Zou type, the entropy
method can be replaced by Lyapunov--Krasovskii functionals. For the
Perron--Volterra constructions here, the main difficulty is to build a
Krasovskii functional whose delayed terms preserve the signed invasion
structure on every face. The persistence and boundary-relay arguments also
have to be replaced by their delay-dynamical analogues.

For reaction--diffusion equations with homogeneous coefficients and
Neumann boundary conditions, spatially constant equilibria preserve the
same algebraic reproduction numbers. Entropy integrals often produce
additional nonpositive diffusion terms, for example terms of the form
\[
-\int_\Omega D_i z_i^*
\frac{|\nabla u_i|^2}{u_i^2}\,dx.
\]
Thus the entropy bracket is not the main obstruction in the simplest
homogeneous Neumann case. The difficulty is the spectral and persistence
theory for the full PDE semiflow, especially with heterogeneous
coefficients or different diffusion rates.

For stochastic finite-population models, deterministic GAS has no literal
analogue, because extinction occurs with probability one in finite state
spaces with absorbing boundary states. The deterministic reproduction and
invasion numbers retain their local threshold role in large-population
limits: near a boundary equilibrium, the missing strain is approximated by
a branching or multitype branching process whose mean growth is governed
by the corresponding transversal Metzler block. A stochastic CEP analogue
would therefore have to be formulated in terms of large-population
metastability, quasi-stationary distributions, fixation probabilities, or
most likely invasion paths, not as global convergence to a deterministic
equilibrium.
\eeR